\documentclass[MSNbibl,number,citesort,seceqn,dvips]{arxbj}
\usepackage[left]{vmkcol}
\usepackage{dcolumn,upgreek,accents}

\aid{0}
\volume{18}
\issue{2}
\pubyear{2012}
\firstpage{644}
\lastpage{677}
\doi{10.3150/10-BEJ349}

\makeatletter
\newcolumntype{d}[1]{D{.}{.}{#1}}
\newproclaim{assumption}{Assumption}
\newtheorem{theorem}{Theorem}
\newtheorem{lemma}{Lemma}
\newtheorem{proposition}{Proposition}

\makeatother

\begin{document}
\begin{frontmatter}

\title{Inference on power law spatial trends}
\runtitle{Inference on power law spatial trends}

\begin{aug}
\author{\fnms{Peter M.} \snm{Robinson}\corref{}\ead[label=e1]{p.m.robinson@lse.ac.uk}}
\runauthor{P.M. Robinson}
\address{Department of Economics, London School of Economics, London
WC2A 2AE, UK.\\
\printead{e1}}
\end{aug}

\received{\smonth{1} \syear{2010}}
\revised{\smonth{12} \syear{2010}}

%
\begin{abstract}
Power law or generalized polynomial regressions with unknown real-valued
exponents and coefficients, and weakly dependent errors, are considered for
observations over time, space or space--time. Consistency and asymptotic
normality of nonlinear least-squares estimates of the parameters are
established. The joint limit distribution is singular, but can be
used as
a basis for inference on either exponents or coefficients. We discuss
issues of implementation, efficiency, potential for improved estimation and
possibilities of extension to more general or alternative trending
models to allow for irregularly spaced data or heteroscedastic errors;
though
it focusses on a particular model to fix ideas, the paper can be viewed as
offering machinery useful in developing inference for a variety of
models in
which power law trends are a component. Indeed, the paper also makes a
contribution that is potentially relevant to many other statistical models:
Our problem is one of many in which consistency of a vector of parameter
estimates (which converge at different rates) cannot be established by the
usual techniques for coping with implicitly-defined extremum estimates, but
requires a more delicate treatment; we present a generic consistency result.
\end{abstract}

%

\end{frontmatter}
%
\section{Introduction}\label{sec1}

Polynomial-in-time regression is one of the longest-established tools of
time series analysis (see Jones~\cite{Jon43}). In much empirical work,
especially when stochastic trends, such as unit roots, are also involved,
only a linear trend is countenanced, or merely a constant intercept. On the
other hand, classical methods can test polynomial order when observations
are equally spaced in time. With independent and identically distributed
(i.i.d.) normal errors, a particularly elegant way of achieving this, with
finite sample validity, results from an orthogonal polynomial representation~-- the covariance
matrix of the least-squares estimate (LSE) is diagonalized,
and contributions to the~$F$ statistic from individual regressors are i.i.d.
(see Section~3.2.2 of Anderson \cite{And71}). Asymptotic theory is valid under
much wider conditions on the errors; indeed from Section~7.4 of Grenander and
Rosenblatt \cite{GreRos84}, the LSE is asymptotically efficient (in the Gauss--Markov
sense) when the (possibly non-Gaussian) errors are covariance stationary
with spectral density bounded and bounded away from zero at zero frequency,
as with short memory processes. Polynomial models have also been
extended to
spatial lattice data (see Section~3.4 of Cressie \cite{Cre93}).

Polynomials are nevertheless restrictive. The Weierstrass theorem justifies
their uniform approximation of any continuous function over a compact
interval, but seems less practically relevant the longer the data set.
Nonparametric smoothing may be unreliable in a series of moderate
length, when
instead richer parametric models than polynomials might be considered. One
class that advantageously nests polynomials, which has received little
theoretical attention, consists of ``generalized polynomial'' or ``power law''
models. With equally spaced time series observations $y_{u}$, $u=1,\ldots,N$,
consider %
\begin{equation}
y_{u}=\sum^{p}_{j=1}\beta_{j}u^{\theta_{j}}+x_{u},
\label{eq1.1}
\end{equation}
where the $\theta_{j}$ and $\beta_{j}$ are real valued and all can be
unknown, $\theta_{j}>-1/2$ for all $j$, and the zero-mean unobservable
process $x_{u}$ is covariance stationary with short memory. For $\theta
_{j}<-1/2$, $\beta_{j}$ would not be estimable (whether $\theta_{j}$ were
known or unknown) because the corresponding signal is drowned by the noise.
For $\theta_{j}=-1/2$, $\beta_{j}$ is estimable but we omit this
possibility because our central limit theorem requires $\theta_{j}$ to lie
in the interior of a~compact set. Polynomials, such as when $\theta_{j}=j-1$
for all $j$, are nested; indeed this is a~hypothesis that might be tested
within~(\ref{eq1.1}).

We consider the nonlinear least-squares estimate (NLSE) of the $\theta
_{j},\beta_{j}$ in (\ref{eq1.1}) and, more generally, of exponents and
coefficients in an extended model defined on a lattice, applying to spatial
and spatio-temporal data, where our provision, for example, for weaker
trends than linear ones and for decaying trends seems practically
useful.
Unlike the LSE when exponents are known, the NLSE cannot be expressed in
closed form and requires numerical optimization.  Correspondingly,
asymptotic theory, with sample size~$N$ increasing, is needed to justify
rules of statistical inference even when errors are Gaussian.  We establish
consistency and asymptotic normality for the NLSE of exponent and
coefficient estimates, achieving also an analogous efficiency bound to that
described above.  As with other implicitly defined estimates, asymptotic
distribution theory makes use (in application of the mean value
theorem) of
an initial consistency proof.  Many such proofs (see Jennrich~\cite{Jen69}, Malinvaud \cite{Mal70})
require regressors to be non-trending,  whence under
suitable additional conditions all parameter estimates are $N^{1/2}$-consistent.
For the NLSE of (\ref{eq1.1}), Wu \cite{Wu81} significantly relaxed this
requirement but nevertheless appears to heavily restrict the diversity of
trends. The discussion after Assumptions A and~A$^\prime$ of Wu \cite{Wu81} indicates
that they reduce in (\ref{eq1.1}) with known $\theta_{j}$ to the assumption
$\max_{j}\theta_{j}<{\frac12}+2\min_{j}\theta_{j}$,
and no weaker requirement suffices in the case of
unknown~$\theta_{j}$.  Example 4 of Wu \cite{Wu81} addressed the latter case
but with $p=1$ only (and for $\theta_{1}\in(-{\frac12},0]$) when the inequality
is trivially satisfied. In general, more elaborate
techniques seem required to establish consistency in (\ref{eq1.1}). Moreover, Wu
\cite{Wu81} established consistency with no rate, whereas we find that a slow
rate of convergence in the $\theta_{j}$ estimates is required before
asymptotic normality is established.  Wu \cite{Wu81} also established asymptotic
normality of the NLSE in a quite general setting, but under the assumption
that all parameter estimates converge at the same rate.  This is not the
case with (\ref{eq1.1}); indeed all rates of $\theta_{j},\beta_{j}$ estimates turn
out to differ. For implicitly defined extremum estimates such variation is
typically associated with difficulty in the initial consistency proof due
to the objective function not converging uniformly to a function that is
uniquely optimized over the whole parameter space.  Consistency proofs here
have tended to be geared to the case at hand (see e.g. Giraitis, Hidalgo and Robinson \cite{GirHidRob01},
Nagaraj and Fuller \cite{NagFul91}, Nielsen \cite{Nie07}, Robinson \cite{Rob08}, Sun and Phillips~\cite{SunPhi03}). Our consistency proof
employs a generic result (presented and proved in Appendix \hyperref[appa]{A} to avoid
interrupting the flow) that seems likely to apply to a quite general class
of estimates (not just the NLSE) of a variety of models.  Our asymptotic
distribution theory of estimates for (\ref{eq1.1}) and its extension presents some
other unusual features.

The following section presents the model, regularity conditions and three
theorems describing asymptotic statistical properties.  The main
details of
their proofs appear in Appendix \hyperref[appb]{B}. These use a series of propositions,
stated and proved in Appendix \hyperref[appc]{C}, and relying in turn also on a series of
lemmas, in Appendix \hyperref[appd]{D}. A Monte Carlo study of finite sample performance
appears in Section \ref{sec3}, while Section \ref{sec4} discusses aspects of the
theoretical results and
their implementation, with possible extensions.

\section{Estimation of spatial lattice regression model}\label{sec2}

Let the integer $d\geq1$ represent the dimension on which data are
observed, $ $where $d=1$ for time series (as in (\ref{eq1.1})) and $d\geq2$ for
spatial or spatio-temporal data. Generalize $u$ to the $d$-dimensional
multi-index $u=( u_{1},u_{2},\ldots,u_{d}) ^{\prime}$.  Denoting $%
\mathbb{Z}_{+}=\{ j\dvt j=0,1,\ldots\} $, generalize~(\ref{eq1.1}) to
%
\begin{equation}\label{eq2.1}
y_{u}=\sum^{d}_{i=1}\sum^{p_{i}}_{j=1}\beta_{ij}u_{i}^{\theta_{ij}}+x_{u}=f(u;\theta
)^{\prime
}\beta+x_{u},\qquad  u\in\mathbb{Z}_{+}^{d},
\end{equation}
where $x_{u}$ is described subsequently and
$\beta=( \beta_{1}^{\prime},\ldots,\beta_{d}^{\prime}) ^{\prime}$,
$\beta_{i}=( \beta_{i1},\ldots,\beta_{ip_{i}}) ^{\prime}$,
$\theta =( \theta_{1}^{\prime},\allowbreak\ldots, \theta_{d}^{\prime})^{\prime}$,
$\theta_{i}= ( \theta_{i1},\ldots,\theta_{ip_{i}}) ^{\prime}$,
$f(u;\theta)=(f_{1}(u_{1};\theta_{1})^{\prime},\ldots,f_{d}(u_{d};\theta_{d})^{\prime})^{\prime}$,
$f_{i}(u_{i};\theta_{i})=( u_{i}^{\theta_{i1}},\ldots,\allowbreak u_{i}^{\theta_{ip_{i}}})^{\prime}$,
for $i=1,\ldots,d$. Defining $p=p_{1}+\cdots+p_{d}$,
the $p\times1$ vectors $\beta$ and $\theta$ are supposed unknown.
Any $f_{i}(u_{i};\theta_{i})$ might be absent from $f(u;\theta)$ when
corresponding $\theta_{i}$ and $\beta_{i}$ are void; we proceed as if
corresponding $p_{i}$ and sums over $j=1,\ldots,p_{i} $ are zero,
avoiding indicator functions to describe such circumstances.

Our consistency proof confines the NLSE of $\theta$ to a compact set.
Prescribe an (arbitrarily small) positive $\delta$, and for each
$i=1,\ldots,d$%
, prescribe $\underaccent\bar{\Delta}_{i}$, $\bar{\Delta}_{i}$ such
that $%
-1/2<\underaccent\bar{\Delta}_{i}<\bar{\Delta}_{i}<\infty$, and define
%
\begin{equation} \label{eq2.2}
\Theta_{i}=\{ h_{1},\ldots,h_{p_{i}}\dvt h_{1}\geq\underaccent\bar{\Delta
}_{i};%
h_{j}-h_{j-1}\geq\delta, j=2,\ldots,p_{i};h_{p_{i}}\leq
\bar{\Delta}_{i}\}
\end{equation}
and $\Theta=$ $\prod_{i=1}^{d}\Theta_{i}\mathsf{.}$ We introduce two
assumptions that imply identifiability of $\theta$ and $\beta$.

\begin{assumption}\label{ass1}
$\theta\in\Theta$.
\end{assumption}

\begin{assumption}\label{ass2}
$\theta_{ij}=0$ for at most one $(i,j) ;$ $\beta_{ij}\neq0$ for all $( i,j)
$.\vadjust{\goodbreak}
\end{assumption}

Assumption \ref{ass1} implies%
%
\begin{equation}
-1/2<\theta_{i1}<\cdots<\theta_{ip_{i}}<\infty, \qquad i=1,\ldots,d.
\label{eq2.3}
\end{equation}
The ordering in (\ref{eq2.3}) is arbitrary, and distinctness of the $\theta_{ij}$
across $j$ along with the first part of Assumption \ref{ass2} identifies $\beta$;
note that $u_{i}^{0}=1$ for all $i$ and that we allow an intercept but do
not require one. The second part of Assumption \ref{ass2} identifies $\theta$.

Given $N=\prod _{i=1}^{d}n_{i}$ observations on $y_{u}$, $u\in
\mathbb{N=N}_{1}\times\cdots\times\mathbb{N}_{d}$  and $\mathbb{N}%
_{i}=(1,\ldots,n_{i})$, define the NLSE of $\beta$, $\theta$ by $(
\hat{\beta},\hat{\theta}) =\arg\min_{b\in\mathbb{R}%
^{p},h\in\Theta}Q(b,h)$, where $Q(b,h)=\sum _{u\in\mathbb
{N}%
}\{ y_{u}-b^{\prime}f(u;h)\} ^{2}$.  Asymptotic theory
requires further assumptions. Let $\mathbb{Z=}\{ j\dvt j=0,\pm
1,\ldots\} $.

\begin{assumption}\label{ass3}
$\mathit{x}_{u}$, $u\in\mathbb{Z}^{d}$,
is covariance stationary with zero mean, and its autocovariance
function, $\gamma_{u}=\operatorname{cov}( x_{t},x_{t+u})$, for the
multi-index $t=( t_{1},\ldots,t_{d}) ^{\prime}$, satisfies
$\sum _{u\in\mathbb{Z}^{d}}\vert\gamma_{u}
\vert <\infty$.
\end{assumption}

Our parameter estimates make no attempt to correct for this possible
nonparametric weak dependence of the $x_{u}$ (permitted also in Assumption
\ref{ass5}), and Cressie \cite{Cre93}, page~25, stresses the importance of mean
function specification relative to error specification.  However, the
NLSE turns out to be not only consistency-robust to spatial correlation but also
asymptotically Gauss--Markov efficient.

The next assumption, of increase with algebraic rate of observations in all
dimensions, is capable of generalization but is employed for
simplicity.

\begin{assumption}\label{ass4}
$n_{i}\sim B_{i}N^{b_{i}}$, $i=1,\ldots,d$,
as $N\rightarrow\infty$, where $B_{i}>0$, $b_{i}>0$,
$i=1,\ldots,d$, $\prod _{i=1}^{d}B_{i}=\sum_{i=1}^{d}b_{i}=1$.
\end{assumption}

Define $\zeta_{ij}=b_{i}\theta_{ij}$ and, with no loss of generality,
identify dimension $i=1$ such that
%
\begin{equation}\label{eq2.4}
\zeta_{11}=\min_{1\leq i\leq d}\{ \zeta_{i1}\} ,
\end{equation}
where, if two or more $i$ satisfy (\ref{eq2.4}), an arbitrary choice is made. Note
that $\zeta_{11}+{\frac12}>0$ is implied by $\theta_{11}+{\frac12}>0$.

\begin{theorem}\label{th1}
Let Assumptions \textup{\ref{ass1}--\ref{ass4}} hold. Then for $j=1,\ldots,p_{i}$, $i=1,\ldots,d$, as $%
N\rightarrow\infty$,%
%
\begin{equation}\label{eq2.5}
\hat{\theta}_{ij}-\theta_{ij}=\mathrm{O}_{p}( N^{\chi-\zeta_{ij}-{1/2}})
\end{equation}
for any $\chi>0$.
\end{theorem}

The proof is in Appendix \hyperref[appb]{B}. As is common with initial consistency
proofs, a
sharp rate (corresponding to $\chi=0$ in (\ref{eq2.5})) is not delivered
(smoothness conditions, in particular, are not exploited). Theorem \ref{th1} is used
in the proof of our central limit theorem (CLT), for which we also need
consistency, with a rate, for $\hat{\beta}$. We state this result without
the proof, which is a relatively straightforward application of Theorem \ref{th1},
techniques used in its proof, Theorem \ref{th3} below and routine
manipulations.

\begin{theorem}\label{th2}
Let Assumptions \ref{ass1}--\ref{ass4} hold. Then, for $j=1,\ldots,p_{i}$, $i=1,\ldots,d$,%
\[
\hat{\beta}_{ij}=\beta_{ij}+\mathrm{O}_{p}( (\log N)N^{\chi-\zeta_{ij}-{1/2}
}),\qquad \mbox{as }N\rightarrow\infty.
\]
\end{theorem}

The \textit{relative} rates for the $\hat{\theta}_{ij}$ and $\hat{\beta}
_{ij} $ in Theorems \ref{th1} and \ref{th2} are matched by relative rates that feature in
our CLT. For this we introduce first

\begin{assumption}\label{ass5}
$x_{u}=\sum _{v\in\mathbb
{Z}^{d}}\xi_{v}\varepsilon_{u-v}$,
$\sum _{v\in\mathbb{Z}^{d}}\vert\xi_{v}\vert<\infty$,
$u\in\mathbb{Z}^{d}$,\vspace*{2pt}
where ${v}$ is the multi-index $v=(v_{1},\ldots,v_{d})^{\prime}$,
$\{ \varepsilon_{u},u\in\mathbb
{Z}^{d}\}$ are independent random variables with zero
mean and
unit variance, $\{ \varepsilon_{u}^{2},u\in\mathbb
{Z}^{d}\} $ are uniformly integrable and $\sum _{v\in\mathbb
{Z}%
^{d}}\xi_{v}\neq0$.
\end{assumption}

Assumption \ref{ass5} implies Assumption \ref{ass3}, and both imply the existence and boundedness
of the spectral density $F(\lambda)=( 2\uppi) ^{-1}
\vert
\sum _{v\in\mathbb{Z}^{d}}\xi_{v}\mathrm{e}^{\mathrm{i}v^{\prime}\lambda
}
\vert^{2}$ of $x_{u}$, where $\mathit{\lambda}$ is the multi-index $%
\mathit{\lambda}=(\mathit{\lambda}_{1},\ldots,\mathit{\lambda
}_{d})^{\prime
}$, while Assumption \ref{ass5} also implies $F(0)>0$.  Stationary invertible
autoregressive moving averages are among time series processes covered by
Assumption \ref{ass5}, as are spatial generalizations of these
(see e.g. Hallin, Lu and Tran \cite{HalLuTra01},
Robinson and Vidal~Sanz \cite{RobVid06},
Tj{\o}stheim \cite{Tjs78,Tjs83},
Yao and Brockwell \cite{YaoBro06}).
Mixing conditions, such as ones employed
in a
spatial context by Gao, Lu and Tj\o stheim \cite{GaoLuTjs06}, Hallin, Lu and Yu
\cite{HalLuYu09}, and
Lu, Lundervold, Tj\o stheim and Yao~\cite{Luetal07}, provide an alternative
route for
establishing a CLT, but are not strictly weaker or
stronger than Assumption \ref{ass5}, which we prefer here because $x_{u}$, unlike
processes considered in the latter references, is involved only linearly.

Let $I_{r}$ be the $r$-rowed identity matrix, $\otimes$ denote the Kronecker
product, and introduce $p\times p$ matrices $D=N^{1/2}\operatorname{diag}\{ n_{1}^{\theta_{11}},\ldots,n_{1}^{\theta
_{1p_{1}}},\ldots,n_{d}^{\theta_{d1}},\ldots,n_{d}^{\theta_{dp_{d}}}
\}$, $L(s)= \operatorname{diag}\{ L_{1}(s_{1}),\ldots,L_{d}(s_{d})\}$,
where $L_{i}(s_{i})=(\log s_{i})I_{p_{i}}$, and $(2p\times2p)$
matrices $D_{+}=I_{2}\otimes D$ and $L_{+}=\operatorname{diag}\{ I_{p},L(n)\}$.
Define $\alpha=(\theta^{\prime},\beta^{\prime})^{\prime}$, $\hat
{\alpha}=(%
\hat{\theta}^{\prime},\hat{\beta}^{\prime})^{\prime}$. Denote
by $\mathfrak{N}_{r}(a,A)$ an $r$-dimensional normal vector with mean vector
$a$ and (possibly singular) covariance matrix~$A$. Appendix \hyperref[appb]{B} defines
the $p\times p$ matrix $\Upsilon$ and $p\times2p$ matrix $B$ and proves:

\begin{theorem}\label{th3}
Let Assumptions \ref{ass1}, \ref{ass2} and \ref{ass5} hold. Then as $N\rightarrow\infty$,
\[
D_{+}L_{+}^{-1}( \hat{\alpha}-\alpha) \rightarrow_{d}%
\mathfrak{N}_{2p}( 0,2\uppi F(0)B^{\prime}\Upsilon^{-1}B) .
\]
\end{theorem}

\section{Finite sample properties}\label{sec3}

A small Monte Carlo study provides some information on finite sample
performance. Issues of concern, given unknown $\theta$, are bias and
variability of the NLSE and accuracy of large sample inference rules
suggested by Theorem \ref{th3}. We employed (\ref{eq2.1}) with $d=2$, $p_{1}=p_{2}=1$,
picking $2$ $( \theta_{1},\theta_{2}) =( \theta
_{11},\theta_{21}) $ combinations -- $(1,1)$, $(0.5,2)$ -- but
throughout took $\Theta_{i1}=[-0.45,4]$, $\beta_{i}=\beta_{i1}=1$,
$i=1,2$. We varied $N$ absolutely and also the relative~$n_{1},n_{2}$, taking
$n_{1},n_{2}=(8,12)$, $(10,10),(11,20),(15,15)$.

\begin{table}
  \caption{$\theta_{1}=1$, $\theta_{2}=1$, $\beta_{1}=1$, $\beta
_{2}=1$, $\sigma^{2}=1$, $x_{u}$ i.i.d.}\label{tab1}
\begin{tabular*}{\textwidth}{@{\extracolsep{\fill}}llld{2.3}lld{2.3}d{2.3}l@{}}
\hline
$n_{1}$ & $n_{2}$ & & \multicolumn{1}{l}{$\hat{\theta}_{1}$} & $\hat{\theta}_{2}$
& $\hat{\beta}_{1}$ & \multicolumn{1}{l}{$\tilde{\beta}_{1}$} & \multicolumn{1}{l}{$\hat{\beta}_{2}$} &
$\tilde{\beta}_{2}$\\
\hline
\phantom{0}$8$ & $12$ & BIAS &0.008&0.007&0.024&0.000&0.017&0.000\\
&& MSE &0.016&0.007&0.080&0.001&0.051&0.000 \\
&& SIZE5 &0.100&0.125&0.151&0.048&0.166&0.055\\
&& SIZE1 &0.044&0.048&0.075&0.010&0.084&0.010\\[3pt]
$10$ & $10$ & BIAS& 0.005&0.009&0.016&-0.001&0.009&0.002\\
&&MSE &0.010&0.009&0.060&0.006&0.063&0.007\\
&&SIZE5&0.132& 0.132&0.180&0.053&0.186&0.051\\
&&SIZE1&0.055&0.050&0.084&0.015&0.090&0.011\\[3pt]
$11$ & $20$ &BIAS& -0.002&0.002&0.016&0.000&-0.007&0.000\\
&&MSE& 0.003&0.001&0.022&0.000&0.010& 0.000\\
&&SIZE5& 0.086&0.104&0.115&0.039&0.120&0.051 \\
&&SIZE1&0.030&0.039&0.051&0.005&0.049&0.012\\[3pt]
$15$ & $15$ & BIAS &0.003&0.002&0.006&0.000&-0.001&0.000\\
&&MSE &0.002&0.002&0.013&0.000&0.013&0.000\\
&&SIZE5 &0.074&0.075&0.108&0.043&0.103&0.039 \\
&& SIZE1&0.024&0.022&0.033&0.010&0.037&0.010\\
\hline
\end{tabular*}
\end{table}

Our  first experiment took the $x_{u}$ to be i.i.d. $\mathfrak
{N}_{1}(0,1) $ variables.  Tables \ref{tab1} and \ref{tab2} report, for the respective
parameter combinations, bias (BIAS), mean squared error (MSE), and empirical
size at 5\% (SIZE5) and 1\% (SIZE1) for the NLSE $\hat{\theta}_{i}$,
$\hat{%
\beta}_{i}$, and also $\tilde{\beta}_{i}$, the LSE of $\beta_{i}$ that
correctly assumes $\theta$, $ $for $i=1,2$, across $1000$ replications.
The sizes were proportions of significant estimates, using normal critical
values scaled by estimated standard deviations which, in the case of
the $\hat{%
\theta}_{i}$, $\hat{\beta}_{i}$, were computed on the basis of Theorem
\ref{th3} with
current parameter estimates replacing true values of $\theta,\beta$,
and $%
2\uppi F(0)$ replaced by the sum of squared residuals divided by $N$ (so the
spatial independence of the $x_{u}$ was treated as known, as it was
also in
the conventional scaling used for the $\tilde{\beta}_{i}$).

The tables reveal a definite inferiority of the NLSE relative to the LSE,
but unsurprisingly, as the LSE is exactly unbiased, more efficient and
yields exact critical regions. Though the NLSE-based tests on $\beta$ are
nearly always over-sized, this phenomenon diminishes with increased
$N$, and
overall the discrepancy between the performances of the two classes of
the $%
\beta$ estimate does not seem very serious. There is also a predominate
over-sizing of the tests on $\theta$, but again this falls as $N$
increases, and, in Table \ref{tab2} in particular, it is often modest. There is a
tendency for the NLSE to over-estimate, but for $\beta$ biases only exceed
2\% of the parameter value when $n_{i}=8$ and $n_{i}=12$. For $\theta$
they never reach 1\%, while overall they mostly fall with increasing
$N$, as
does the MSE.  In Table~\ref{tab2}, the results are not in line with what the rates
in Theorem \ref{th3} suggest, because the fall in MSE is greater for $\hat
{\theta}_{2}$ and $\hat{\beta}_{2}$ than for $\hat{\theta}_{1}$ and $\hat{\beta
}_{1}$, despite the fact that $\theta_{1}=2$ and $\theta_{2}={\frac12}$.
Nevertheless, it is not clear to what extent one would expect
asymptotic theory to predict comparisons at this level of refinement in such
sample sizes. Note that the Monte Carlo results are also difficult to
judge relative to the theory because the various~$n_{i}$ did not result from
fixing the $b_{i}$ and $B_{i}$ and then increasing $n$, but were chosen with
a view to representing some variability in $n$, and relative to $n_{1}$
and $n_{2}$. In addition, the convergence rates of $\hat{\theta}_{i}$
and $\hat{\beta}_{i}$ do not only depend on $n_{i}$, but on the overall $n$.
Other results are more closely in line with the asymptotic theory.
This is the case in Table~\ref{tab1} where, with $\theta_{1}=\theta_{2}=1$, the
above MSE ratios are sometimes greater for $\hat{\theta}_{2}$
and/or~$\hat{\beta}_{2}$ and sometimes less. It is also the case in Table \ref{tab2} for the
LSE $\tilde{\beta}_{i}$, as elsewhere, that comparisons are sometimes
difficult as a number of MSEs are zero to $3$, and even to $4$ (unreported
here), decimal places.

\begin{table}
\caption{$\theta_{1}=2$, $\theta_{2}=1/2$, $\beta_{1}=1$, $\beta
_{2}=1$, $\sigma^{2}=1$, $x_{u}$ i.i.d.}\label{tab2}
\begin{tabular*}{\textwidth}{@{\extracolsep{\fill}}lllld{2.3}d{2.3}d{2.3}d{2.3}d{2.3}@{}}
\hline
\multicolumn{1}{@{}l}{$n_{1}$} &
\multicolumn{1}{l}{$n_{2}$} & & \multicolumn{1}{l}{$\hat{\theta}_{1}$} &
\multicolumn{1}{l}{$\hat{\theta}_{2}$}
&
\multicolumn{1}{l}{$\hat{\beta}_{1}$} &
\multicolumn{1}{l}{$\tilde{\beta}_{1}$} &
\multicolumn{1}{l}{$\hat{\beta}_{2}$} &
\multicolumn{1}{l@{}}{$\tilde{\beta}_{2}$}\\
\hline
\phantom{0}$8$ & $12$ & BIAS& 0.008&0.001&0.024&0.003&-0.002&-0.000\\
&&MSE& 0.014&0.001&0.071&0.005&0.001&0.000\\
&&SIZE5& 0.063&0.060&0.087&0.077&0.053&0.090\\
&&SIZE1&0.028&0.012&0.038&0.029&0.014&0.034\\[3pt]
$10$ & $10$ &BIAS &0.008&0.000&0.020&0.004&0.000&-0.000\\
&&MSE &0.013&0.003&0.074&0.004&0.001&0.000\\
&&SIZE5&0.069&0.057 &0.101&0.058&0.065&0.039\\
&&SIZE1&0.033&0.013&0.047&0.015&0.017&0.009\\[3pt]
$11$ & $20$ & BIAS &0.005&-0.000&-0.001&-0.002&0.000&0.000\\
&&MSE &0.005&0.000&0.028&0.002&0.000 &0.000\\
&&SIZE5 &0.052&0.054&0.069&0.030&0.059&0.041\\
&& SIZE1&0.017&0.012&0.017&0.012&0.011&0.006\\[3pt]
$15$ & $15$ & BIAS&0.002 &0.001&0.004&0.001&0.004&0.000\\
&&MSE &0.004&0.001&0.025&0.001&0.000&0.000 \\
&&SIZE5& 0.058&0.044&0.070&0.081&0.043 &0.055\\
&& SIZE1&0.018&0.011&0.019&0.019&0.010&0.020\\
\hline
\end{tabular*}
\end{table}

Next we considered the effect of dependence, employing three different
models for $x_{u}$, again with $d=2$. All models entailed weak dependence,
with varying spans, but in the first dependence was negative, so that the
spectral density at zero was small, whereas in the other two it was
positive, producing a peaked spectral density. In the following, $%
\varepsilon_{u}\sim$ i.i.d. $\mathfrak{N}_{1}( 0,1) $.

\begin{enumerate}
\item Multiple direction $\operatorname{MA}(1)$:
%
\begin{equation}\label{eq3.1}
x_{u}=\varepsilon_{u}-0.12
\mathop{\sum _{j=-1}^{1}\sum _{k=-1}^{1}}_{(j,k)\neq0}
\varepsilon_{u_{1}+j,u_{2}+k},\qquad  u_{i}=1,\ldots,n_{i}, i=1,2.
\end{equation}

\item Multilateral $\operatorname{MA}(4)$, no interactions:
%
\begin{equation}\label{eq3.2}
x_{u}=\varepsilon_{u}+\mathop{\sum_{j=-4}}_{j\neq 0}^{4}a_{\vert j\vert}( \varepsilon_{u_{1}+j,u_{2}}+\varepsilon
_{u_{1},u_{2}+j}) ,\qquad  u_{i}=1,\ldots,n_{i}, i=1,2
\end{equation}
for $a_{1}=0.14$, $a_{2}=0.12$, $a_{3}=0.1$, $a_{4}=0.08$.

\item Bilateral $\operatorname{MA}(9)$, on diagonal:
%
\begin{equation}\label{eq3.3}
x_{u}=\varepsilon_{u}+\mathop{\sum_{j=-9}}_{j\neq0}^{9}(0.95)^{\vert j\vert}\varepsilon
_{u_{1}+j,u_{2}+j},\qquad u_{i}=1,\ldots,n_{i}, i=1,2.\vspace*{-3pt}
\end{equation}
\end{enumerate}

\begin{table}
\caption{$\theta_{1}=1$, $\theta_{2}=1$, $\beta_{1}=1$, $\beta
_{2}=1$, $\sigma
^{2}=1$, $x_{u} =$ (\protect\ref{eq3.1})}\label{tab3}
\begin{tabular*}{\textwidth}{@{\extracolsep{\fill}}llld{1.3}d{2.4}d{2.3}d{2.3}ld{2.3}@{}}
\hline
\multicolumn{1}{@{}l}{$n_{1}$} &
\multicolumn{1}{l}{$n_{2}$} & & \multicolumn{1}{l}{$\hat{\theta}_{1}$} &
\multicolumn{1}{l}{$\hat{\theta}_{2}$}
&
\multicolumn{1}{l}{$\hat{\beta}_{1}$} &
\multicolumn{1}{l}{$\tilde{\beta}_{1}$} &
\multicolumn{1}{l}{$\hat{\beta}_{2}$} &
\multicolumn{1}{l@{}}{$\tilde{\beta}_{2}$}\\
\hline
\phantom{0}$8$ & $12$ & BIAS &0.005&0.003&0.006&0.000&0.002&-0.000\\
&& MSE&0.006&0.003&0.031&0.000&0.021&0.000\\
$10$ & $10$ & BIAS&0.005 &0.001&0.001&0.000&0.008&-0.000\\
&& MSE&0.003&0.003&0.023&0.000&0.023&0.000\\[3pt]
$11$ & $20$ & BIAS&0.001& 0.0001&0.001&-0.000&0.001&0.000\\
&& MSE&0.001&0.000&0.006&0.000&0.003&0.000\\[3pt]
$15$ & $15$ & BIAS &0.002&-0.001&-0.003&-0.000&0.005&0.000\\
&& MSE&0.000&0.000&0.004&0.000&0.004&0.000\\
\hline
\end{tabular*}
\vspace*{-3pt}
\end{table}

\begin{table}[b]
\vspace*{-3pt}
\caption{$\theta_{1}=2$, $\theta_{2}=1/2$, $\beta_{1}=1$, $\beta
_{2}=1$, $\sigma^{2}=1$, $x_{u} =$ (\protect\ref{eq3.1})}\label{tab4}
\begin{tabular*}{\textwidth}{@{\extracolsep{\fill}}llld{2.3}lld{2.3}d{2.3}d{2.3}@{}}
\hline
\multicolumn{1}{@{}l}{$n_{1}$} &
\multicolumn{1}{l}{$n_{2}$} & & \multicolumn{1}{l}{$\hat{\theta}_{1}$} &
\multicolumn{1}{l}{$\hat{\theta}_{2}$}
&
\multicolumn{1}{l}{$\hat{\beta}_{1}$} &
\multicolumn{1}{l}{$\tilde{\beta}_{1}$} &
\multicolumn{1}{l}{$\hat{\beta}_{2}$} &
\multicolumn{1}{l@{}}{$\tilde{\beta}_{2}$}\\
\hline
\phantom{0}$8$ & $12$ & BIAS& 0.003&0.000&0.003&-0.000&-0.000&0.000\\
&& MSE&0.003&0.000&0.017&0.001&0.000&0.000\\[3pt]
$10$ & $10$ & BIAS &-0.003&0.000&0.014&-0.001&-0.001&0.000\\
&& MSE&0.003&0.000&0.018&0.001&0.000&0.000\\[3pt]
$11$ & $20$ & BIAS &-0.000&0.000&0.003&0.000&-0.000& -0.000\\
&& MSE&0.001&0.000&0.004&0.000&0.000&0.000\\[3pt]
$15$ & $15$ & BIAS &-0.001&0.000&0.005&0.001&-0.000&-0.000\\
&& MSE&0.000&0.000&0.004&0.000&0.000&0.000\\
\hline
\end{tabular*}
\end{table}

For the same parameter values as before, bias and MSE of the LSE and NLSE
are presented in Tables \ref{tab3}--\ref{tab8}, with Tables \ref{tab3} and \ref{tab4} referring to (\ref{eq3.1}), Tables
\ref{tab5} and \ref{tab6} to (\ref{eq3.2}), and Tables \ref{tab7} and~\ref{tab8} to (\ref{eq3.3}). As before the LSE $\tilde
{\beta}_{1},\tilde{\beta}_{2}$ are exactly unbiased, as the Monte Carlo
results tend to illustrate. However, perhaps surprisingly, the dependent
model (\ref{eq3.3}) produces some very large biases in the NLSE $\hat{\beta}_{1}$,
though not so much in $\hat{\beta}_{2},\hat{\theta}_{1},\hat{\theta}_{2}$.
For the other dependence models the NLSE biases are not necessarily
greater than under independence. The MSE magnitudes are not directly
comparable to those of Tables \ref{tab1} and \ref{tab2}, because scales were not calibrated,
but a similar overall picture emerges: the NLSE of $\beta$ often has much
greater MSE than the LSE, but this falls with increasing $N$, as does that
of the NLSE of $\theta$.  In Tables \ref{tab4}, \ref{tab6} and \ref{tab8}, where $\theta
_{1}=2$, $\theta_{2}={\frac12}$,
the same somewhat surprising feature\vadjust{\goodbreak} as noted in Table \ref{tab2} appears,
with $\hat{\theta}_{1}$ and $\hat{\beta}_{1}$ improving less than $\hat{\theta
}_{2} $ and $\hat{\beta}_{2}$ with increasing $n$, and the only additional
point to add to our previous discussion is that convergence is often
expected to be slowed by dependence.

\begin{table}
\caption{$\theta_{1}=1$, $\theta_{2}=1$, $\beta_{1}=1$, $\beta
_{2}=1$, $\sigma
^{2}=1$, $x_{u} = $ (\protect\ref{eq3.2})}\label{tab5}
\begin{tabular*}{\textwidth}{@{\extracolsep{\fill}}lllllld{2.3}ll@{}}
\hline
\multicolumn{1}{@{}l}{$n_{1}$} &
\multicolumn{1}{l}{$n_{2}$} & & \multicolumn{1}{l}{$\hat{\theta}_{1}$} &
\multicolumn{1}{l}{$\hat{\theta}_{2}$}
&
\multicolumn{1}{l}{$\hat{\beta}_{1}$} &
\multicolumn{1}{l}{$\tilde{\beta}_{1}$} &
\multicolumn{1}{l}{$\hat{\beta}_{2}$} &
\multicolumn{1}{l@{}}{$\tilde{\beta}_{2}$}\\
\hline
\phantom{0}$8$ & $12$ & BIAS &0.032&0.020&0.053&-0.001&0.035&0.000\\
&& MSE&0.050&0.026&0.249&0.004&0.169&0.002\\[3pt]
$10$ & $10$ & BIAS &0.029&0.020&0.017&-0.005&0.047&0.003\\
&& MSE&0.031&0.031&0.181&0.003&0.177&0.003\\[3pt]
$11$ & $20$ & BIAS &0.010&0.003&0.017&-0.001&0.015&0.001\\
&&MSE&0.013&0.004&0.091&0.001&0.045&0.000\\[3pt]
$15$ & $15$ & BIAS &0.008&0.007&0.006&-0.001&0.014&0.000\\
&& MSE&0.007&0.008&0.059&0.000&0.060&0.001\\
\hline
\end{tabular*}\vspace*{-4pt}
\end{table}

\begin{table}[b]\vspace*{-4pt}
\caption{$\theta_{1}=2$, $\theta_{2}=1/2$, $\beta_{1}=1$, $\beta
_{2}=1$, $\sigma^{2}=1$, $x_{u} =$ (\protect\ref{eq3.2})}\label{tab6}
\begin{tabular*}{\textwidth}{@{\extracolsep{\fill}}lllld{2.3}ld{2.3}d{2.3}d{2.3}@{}}
\hline
\multicolumn{1}{@{}l}{$n_{1}$} &
\multicolumn{1}{l}{$n_{2}$} & & \multicolumn{1}{l}{$\hat{\theta}_{1}$} &
\multicolumn{1}{l}{$\hat{\theta}_{2}$}
&
\multicolumn{1}{l}{$\hat{\beta}_{1}$} &
\multicolumn{1}{l}{$\tilde{\beta}_{1}$} &
\multicolumn{1}{l}{$\hat{\beta}_{2}$} &
\multicolumn{1}{l@{}}{$\tilde{\beta}_{2}$}\\
\hline
\phantom{0}$8$ & $12$ & BIAS &0.064&0.001&0.048&0.005&-0.001&-0.000\\
&& MSE&0.115&0.000&0.272&0.024&0.003&0.000\\[3pt]
$10$ & $10$ & BIAS &0.067&-0.001&0.023&-0.002&0.005&0.000\\
&&MSE&0.111&0.001&0.267&0.019&0.005&0.000\\[3pt]
$11$ & $20$ & BIAS &0.019&0.000&0.035&0.000&-0.001&0.000\\
&& MSE&0.027&0.000&0.151&0.009&0.000&0.000\\[3pt]
$15$ & $15$  &BIAS &0.008&0.000&0.046&-0.002&-0.001&0.000\\
&& MSE&0.020&0.000&0.143&0.007&0.001&0.000\\
\hline
\end{tabular*}
\end{table}

%

\section{Final comments}\label{sec4}
\begin{longlist}[1.]
\item[1.] For known $\theta$, long-established techniques (see \cite{And71}, Section~2.6) give $D( \hat{\beta}(\theta)-\beta)
\rightarrow_{d}\mathfrak{N}_{p}( 0,2\uppi F(0)\Phi^{-1}) $
(where $\Phi$ is defined near the start of Appendix \hyperref[appb]{B} below), so ignorance
of~$\theta$ incurs not only efficiency loss, but slightly slower
convergence. Theorem \ref{th3} also implies a singularity in the limit distribution,
whose covariance matrix has rank $p$ only.\vspace*{1pt} This is due to bias in $\hat
{\beta}$, which on expansion is seen to have a term linear in $\hat{\theta
}-\theta
$ that dominates the contribution from $\sum_{u\in\mathbb
{N}}f(u;\theta
)x_{u}$. Nevertheless, Theorem \ref{th3} does provide separate inference on
$\beta$
(moreover, one can conduct joint inference that does not cover both
$\theta
_{ij}$ and $\beta_{ij}$ for any $(i,j)$), though, given Assumption \ref{ass1}, we
cannot test zero restrictions on~$\beta$. In our setting, $\beta$ may
be of
less initial interest than $\theta$, and Theorem \ref{th3} allows inference on~$\theta$
with $\hat{\theta}$ converging slightly faster than $\hat
{\beta
}$, and at what appears to be the optimal rate for this problem.

\item[2.] If independence of the $x_{u}$ is not assumed, the limiting covariance
matrix in Theorem~\ref{th3} can be consistently estimated (under additional
conditions) by replacing $F(0)$ by a~parametric or smoothed nonparametric
estimate based on NLSE residuals.

\begin{table}
\caption{$\theta_{1}=1$, $\theta_{2}=1$, $\beta_{1}=1$, $\beta
_{2}=1$, $\sigma
^{2}=1$, $x_{u} = $ (\protect\ref{eq3.3})}\label{tab7}
\begin{tabular*}{\textwidth}{@{\extracolsep{\fill}}lllllld{2.4}ld{2.3}@{}}
\hline
\multicolumn{1}{@{}l}{$n_{1}$} &
\multicolumn{1}{l}{$n_{2}$} & & \multicolumn{1}{l}{$\hat{\theta}_{1}$} &
\multicolumn{1}{l}{$\hat{\theta}_{2}$}
&
\multicolumn{1}{l}{$\hat{\beta}_{1}$} &
\multicolumn{1}{l}{$\tilde{\beta}_{1}$} &
\multicolumn{1}{l}{$\hat{\beta}_{2}$} &
\multicolumn{1}{l@{}}{$\tilde{\beta}_{2}$}\\
\hline
\phantom{0}$8$ & $12$ & BIAS &0.074&0.096&0.154&0.008&0.091&-0.004\\
&& MSE&0.129&0.157&0.738&0.048&0.549&0.024\\[3pt]
$10$ & $10$ & BIAS &0.041&0.069&0.105&-0.008&0.050&0.008\\
&& MSE&0.080&0.097&0.455&0.033&0.371&0.032\\[3pt]
$11$ & $20$ & BIAS &0.016&0.036&0.134&0.0010&0.017&-0.000\\
&& MSE&0.043&0.032&0.462&0.014&0.232&0.005\\[3pt]
$15$ & $15$ & BIAS &0.013&0.024&0.061&-0.003&0.028&0.002\\
&& MSE&0.026&0.026&0.214&0.009&0.182&0.009\\
\hline
\end{tabular*}
\end{table}

\begin{table}[b]
\caption{$\theta_{1}=2$, $\theta_{2}=1/2$, $\beta_{1}=1$, $\beta
_{2}=1$, $\sigma^{2}=1$, $x_{u} =$ (\protect\ref{eq3.3})}\label{tab8}
\begin{tabular*}{\textwidth}{@{\extracolsep{\fill}}llld{2.3}d{2.3}d{2.3}d{2.3}ld{2.3}@{}}
\hline
\multicolumn{1}{@{}l}{$n_{1}$} &
\multicolumn{1}{l}{$n_{2}$} & & \multicolumn{1}{l}{$\hat{\theta}_{1}$} &
\multicolumn{1}{l}{$\hat{\theta}_{2}$}
&
\multicolumn{1}{l}{$\hat{\beta}_{1}$} &
\multicolumn{1}{l}{$\tilde{\beta}_{1}$} &
\multicolumn{1}{l}{$\hat{\beta}_{2}$} &
\multicolumn{1}{l@{}}{$\tilde{\beta}_{2}$}\\
\hline
\phantom{0}$8$ & $12$ & BIAS &0.063&-0.000&0.100&0.014&0.009&-0.000\\
&& MSE&0.518&0.003&1.217&0.291&0.019&0.000\\[3pt]
$10$ & $10$ &BIAS &0.098&-0.000&0.118&0.009&0.008&-0.000\\
&& MSE&0.512&0.003&0.912&0.222&0.016&0.000\\[3pt]
$11$ & $20$ & BIAS &-0.037&-0.002&-0.007&-0.001&0.008&0.000\\
&& MSE&0.275&0.000&1.059&0.128&0.004&0.000\\[3pt]
$15$ & $15$ & BIAS &0.054&0.000&0.128&-0.001&0.001 &0.000\\
&& MSE&0.226&0.000&0.616&0.086&0.003&0.000\\
\hline
\end{tabular*}
\end{table}

\item[3.] The form of the limiting covariance matrix in Theorem \ref{th3}, with dependence
simply reflected in the scale factor $2\uppi F(0)$, suggests that a
generalized NLSE, which corrects parametrically or nonparametrically for
correlation in $x_{u}$, affords no efficiency improvement (cf. Section~7.4 of
Grenander and Rosenblatt \cite{GreRos84}).

\item[4.] On the other hand, our estimates are not Fisher efficient for
non-Gaussian $x_{u}$.  Departures from Gaussianity might be detected by,
for example, nonparametric probability density estimation based on
NLSE residuals; Hallin, Lu and Tran \cite{HalLuTra01} studied density estimation for linear
lattice processes.  More efficient parameter estimates could be
obtained by
$M$-estimation using a correctly parameterized $\varepsilon_{u}$
distribution, or adapting semi-parametrically to a nonparametric one, in
either case employing parametric $\{ \xi_{v}\} $ or
approximating them via a long autoregression.  The extra proof details
would be far from trivial, but convergence rates should be unaffected, with
the limiting covariance matrix of Theorem \ref{th3} simply shrunk by a scalar
factor.

\item[5.] Another extension allows long or negative memory, in $x_{u}$,
bearing in
mind results of Yajima \cite{Yaj88} for (\ref{eq1.1}) with known integer $\theta_{i}$,
and Yajima and Matsuda \cite{YajMat}; this would affect all convergence
rates by
the same scalar factor, the efficiency property in Comment~3 would be
lost, and
negative $\theta_{ij}$, and corresponding $\beta_{ij}$ may not be
estimable.

\item[6.] In an alternative formulation to (\ref{eq1.1}), $u^{\theta_{j}}$ is
replaced by
$( u/N) ^{\theta_{j}}$, confining the regression to the unit
interval, and (\ref{eq2.1}) can be analogously modified.  Consistency is then
much easier to prove, all exponent estimates being $\sqrt
{N}$-consistent.  A similar device is employed in fixed-design nonparametric regression, but
unlike there it is not essential in order to achieve consistency in our
parametric setting, where we find it aesthetically unattractive given
that $%
x_{u}$ is defined on an increasing domain.

\item[7.] The results are straightforwardly extended to allow some $\theta_{ij}$
in (\ref{eq2.1}) to be known; for example, to specify an intercept by $\theta
_{11}=0$, though the norming factor and limit covariance matrix in Theorem \ref{th3} are
affected.

\item[8.] Our notation suggests constant spacing between observations across
all $d$
dimensions, but allowing the interval of observation to vary with dimension
affects each $\beta_{ij}$ by a factor depending also on the
corresponding $\theta_{ij}$, but not the $\theta_{ij}$ themselves.

\item[9.] Irregular spacing of observations, either due to missing data from an
otherwise regular lattice, or with observations occurring anywhere on
$\mathbb{R}^{d}$, can also be considered.  In both of these settings
asymptotic theory requires a degree of regularity in the observation
locations, ruling out situations where observations become too sparse, for
example.  Given this, the extension is relatively simple with
independent $%
x_{u}$.  Under dependence, asymptotic variance formulae will be complicated
by the irregular spacing and the efficiency property of Comment 3 will
be lost. In addition, different kinds of assumptions from ours on the
errors $x_{u} $ may be needed.  In the case of missing data from an otherwise
regular lattice, our Assumptions \ref{ass3} (for consistency) and \ref{ass5} (for asymptotic
normality) should still suffice.  But for observations anywhere on
$\mathbb{R}^{d}$ it would be appropriate to consider an underlying continuous
process. Then, for consistency, a suitable ergodicity property
would be
needed, whereas for asymptotic normality leading possibilities that can
entail weak dependence analogous to that of Assumption \ref{ass5} include suitable
linear functionals of Brownian motion and mixing conditions.

\item[10.] A Bayesian treatment would be worthwhile, with suitable priors placed
on the exponents and possibly also the coefficients.

\item[11.] When $d\geq2$ a more realistic model than (\ref{eq2.1}) might allow interaction
terms, that is, products of powers of $u_{i}$ and $u_{k}$, $i\neq k$.  Our
proof methods are extendable, but from a~practical perspective the
curse of
dimensionality threatens and the issue of parsimonious specification,
already posed by (\ref{eq2.1}), becomes more pressing. A penalized procedure could
be used.

\item[12.] Modified model classes might provide an alternative route to parsimony;
for example, one might take $p_{i}=1$ with $\beta_{i1}u_{i1}^{\theta_{i1}}$
replaced by $\beta_{i1}( u_{i1}+\phi_{i1}) ^{\theta
_{i1}}$ for
known or unknown $\phi_{i1}$ (cf. Example 3 of Wu \cite{Wu81}). Trigonometric
factors might also be incorporated (cf. Section~7.5 of Grenander and Rosenblatt
\cite{GreRos84}).

\item[13.]  For alternative classes of trending model (for example, involving
wavelets), asymptotic estimation theory might be handled by similar
techniques.

\item[14.] An alternative practically relevant modelling of the $x_{u}$
treats them
as heteroscedastic but possibly independent. Broadly similar
proof\vadjust{\goodbreak}
techniques would provide corresponding results to ours, but the NLSE is
less efficient than a suitably weighted estimate.\looseness=1

\item[15.] Though we have focussed on (\ref{eq1.1}) and (\ref{eq2.1}) to fix ideas, our
methods and
theory can be developed to cover models that incorporate power law trends
along with other explanatory variables, both stochastic and non-stochastic,
such as extensions of the nonparametric and semiparametric spatial
regressions considered by Gao, Lu and Tj\o st\-heim~\cite{GaoLuTjs06} and Lu,
Lundervold, Tj\o stheim and Yao \cite{Luetal07}, and so the paper can be viewed as introducing
machinery relevant to a wide variety of settings.\looseness=-1\vspace*{-3pt}
\end{longlist}

\begin{appendix}

\section{Generic consistency theorem}\label{appa}\vspace*{-3pt}

We present a consistency theorem for a general, implicitly defined extremum
estimate under unprimitive conditions that will be checked in the paper's
setting and seem capable of checking in a number of others. As this
appendix is self-contained, there seems no risk of confusion in employing
notations that are similar to those elsewhere in the paper but can have
slightly different meanings. We estimate the $p\times1$ vector
parameter $\theta$, with elements~$\theta_{i}$, $i=1,\ldots,p$, by $\hat{\theta
}=\arg
\min_{h\in\Theta}R(h)$, where $R(h)\dvtx\mathbb{R}^{p}\rightarrow\mathbb{R}$
depends on sample size $N$ and $\Theta$ $\subset\mathbb{R}^{p}$ is a fixed
compact set.  For positive scalars $C_{iw}$, $i=1,\ldots,p$, $w=1,2,\ldots,$
depending on $N$ and such that $C_{iw}\leq C_{i,w+1}$, $i=1,\ldots,p$,
define $C_{w}=(C_{1w},\ldots,C_{pw})^{\prime}$, and%
%
\begin{eqnarray}\label{eqA.1}
\mathcal{N}_{i}(C_{iw}) &=&\{ h_{i}\dvt \vert h_{i}-\theta
_{i}\vert<C_{iw}\} ,\qquad \mathcal{N}(C_{w})=\prod_{i=1}^{p}\mathcal{N}_{i}(C_{iw}),
\nonumber
\\[-10pt]
\\[-10pt]
\nonumber
\bar{\mathcal{N}}(C_{w}) &=& \Theta\setminus\mathcal
{N}(C_{w}),\qquad \mathcal{S}_{w}=\bar{\mathcal{N}}(C_{w})\cap\mathcal{N}(C_{w+1}).
\end{eqnarray}

\renewcommand{\thetheorem}{\Alph{theorem}}
\setcounter{theorem}{0}
\begin{theorem}\label{tha}
Assume:
\begin{longlist}[(ii)]
\item[(i)] $\Theta\subset\mathcal{N}(C_{W+1})$ for a finite
integer $W$ and $N$ sufficiently large;

\item[(ii)] There exist positive $s_{1},\ldots,s_{W}$ and
$U(h)$, $V(h)$ such that $R(h)=R(\theta)+U(h)+V(h)$ and
$s_{1}<\cdots<s_{W}$, and as
$N\rightarrow\infty$, $s_{1}\rightarrow\infty$ and
%
\begin{eqnarray} \label{eqA.2}
P\biggl( \inf_{h\in\mathcal{S}_{w}}\frac{U(h)}{s_{w}}>\eta
\biggr) &\rightarrow&1,\qquad  \mbox{some }\eta>0,\\[-2pt]
\label{eqA.3}
 \sup_{h\in\mathcal{S}_{w}}\frac{\vert V(h)
\vert}{s_{w}}&=&\mathrm{o}_{p}(1).
\end{eqnarray}
Then
\[
\hat{\theta}=\theta+\mathcal{O}_{p}(C_{1}), \qquad \mbox{as }N\rightarrow
\infty,
\]
where $\mathcal{O}_{p}(C_{1})$ is a $p\times1$
vector with $i$th element $\mathrm{O}_{p}(C_{i1})$.
\end{longlist}
\end{theorem}

\begin{pf}
We show that $P( \hat{\theta}\in\bar{\mathcal{N}}(
C_{1})) \rightarrow0$ as $N\rightarrow\infty$. By a standard kind
of argument
\[
P\bigl( \hat{\theta}\in\bar{\mathcal{N}}(C_{1}) \bigr)\leq P\Bigl(
\inf_{h\in\bar{\mathcal{N}}(C_{1})}\{ R(h)-R(\theta)\} \leq0\Bigr) .\vadjust{\goodbreak}
\]
Under (i),
$\bar{\mathcal{N}}(C_{1})\subset\bar{\mathcal{N}}(C_{1})\cap
\mathcal{N}(C_{W+1})=\bigcup_{w=1}^{W}\mathcal{S}_{w}$.  Thus the last
probability is bounded by
\[
\sum_{w=1}^{W}P\biggl( \inf_{h\in\mathcal{S}_{w}}\biggl\{\frac{R(h)-R(\theta)}{s_{w}}\biggr\}\leq0\biggr) \leq
\sum_{w=1}^{W}P\biggl( \sup_{h\in\mathcal{S}_{w}}\frac{\vert V(h)\vert}{s_{w}}\geq\inf_{h\in
\mathcal{S}_{w}}\frac{U(h)}{s_{w}}\biggr) ,
\]
which is bounded by%
%
\begin{equation}\label{eqA.4}
\sum_{w=1}^{W}\biggl\{ P\biggl( \sup_{h\in \mathcal{S}_{w}}\frac{\vert V(h)\vert}{s_{w}}>\eta\biggr)
+P\biggl( \inf_{h\in\mathcal{S}_{w}}\frac{U(h)}{s_{w}}\leq
\eta
\biggr) \biggr\} ,
\end{equation}
which tends to zero on applying (\ref{eqA.2}) and (\ref{eqA.3}).
\end{pf}

Three comments are relevant. (1)  In the setting of the rest of the paper,
$U$ can be chosen nonstochastic but this is not possible in the context of
such stochastic trends as unit roots, where the more general (\ref{eqA.2}) is
useful. (2)  An almost sure convergence version of Theorem \ref{tha} is possible
under suitably strengthened versions of (\ref{eqA.2}) and (\ref{eqA.3}).\vspace*{1pt}
(3)  By
comparison with our decomposition of $\bar{\mathcal{N}}(C_{1})$ into $%
\mathcal{S}_{1},\ldots,\mathcal{S}_{W}$, van de Geer \cite{van00} (see pages 69,
70) employed a ``peeling device'' to obtain an exponential inequality
for $\operatorname{sup}
_{g\in\mathcal{G}}\{ \vert Z_{N}(g)\vert/\tau(g) \}$,
where $Z_{N}(g)$ is a stochastic process, $\tau
(g) $ is a non-negative function and the set $\mathcal{G}$ is ``peeled
off'' as $\bigcup _{j=1}^{J}\mathcal{G}_{j}$, where $\mathcal
{G}%
_{j}=\{ g\in\mathcal{G}\dvt m_{j-1}\leq\tau(g)<m_{j}\} $,
for an
increasing sequence $\{ m_{j}\}, $ and $J$ need not be
finite.  Thus $\operatorname{sup}_{g\in\mathcal{G}_{j}}\{ \vert Z_{N}(g)
\vert/\tau( g) \} \leq\{ \operatorname{sup}_{g\in\mathcal
{G}, \tau(g)<m_{j}}\vert Z_{N}(g)\vert\} /m_{j-1}$ and
only the supremum of the numerator of the original statistic need be
approximated.  There is no denominator there like $\tau( g) $ in our
problem, and our decomposition of $\bar{\mathcal{N}}(C_{1})$ is
designed to
suitably balance $U(h)$ and $V(h)$ on each~$\mathcal{S}_{w}$ to enable
choices of the $s_{w}$ that make all $W$ summands in (\ref{eqA.4})
small.\vspace*{-3pt}

\section{Definitions and proofs of theorems}\label{appb}\vspace*{-3pt}

To define $\Upsilon$, introduce first, for $i=1,\ldots,d$, the
$p_{i}\times1$
vector $\phi_{i}(g_{i})$ with $j$th element $( g_{ij}+1) ^{-1}$
and the $p_{i}\times p_{i}$ matrix $\Phi_{i}(g_{i},h_{i})$ with $(j,k)$th
element $( g_{ij}+h_{ik}+1) ^{-1}$ for $g_{i}=(
g_{i1},\ldots,g_{ip_{i}}) ^{\prime}$, $h_{i}=(
h_{i1},\ldots,h_{ip_{i}}) ^{\prime}$, where $g_{ij}$, $h_{ij}>-1/2$ for
all $i,j$. For $g=(g_{1}^{\prime},\ldots,g_{d}^{\prime})^{\prime}$, $%
h=(h_{1}^{\prime},\ldots,h_{d}^{\prime})^{\prime}$, introduce the
$p\times p$
matrix $\Phi(g,h)$ with $(i,j)$th $p_{i}\times p_{j}$ block $\Phi
_{i}(g_{i},h_{i})$ when $i=j$ and $\phi_{i}(g_{i})\phi
_{j}(h_{j})^{\prime
} $ when $i\neq j$.  Denote $\Phi=\Phi(\theta,\theta)$.    Writing $%
\phi_{i}=\phi_{i}(\theta_{i})$, $\Phi_{i}=\Phi_{i}(\theta
_{i},\theta
_{i})$, define $p\times p$ matrices $\Phi_{+}$, $\Phi_{++}$ with $(i,j)$th
$p_{i}\times p_{j}$ block $\Phi_{i}\circ\Phi_{i}$, $2\Phi_{i}\circ
\Phi
_{i}\circ\Phi_{i}$ when $i=j$ and $\phi_{i}( \phi_{j}\circ\phi
_{j}) ^{\prime}$, $( \phi_{i}\circ\phi_{i}) (
\phi
_{j}\circ\phi_{j}) ^{\prime}$ when $i\neq j$, where ``$\circ$''
denotes the Hadamard product.  Put $\Upsilon=\Phi_{++}-\Phi
_{+}^{\prime
}\Phi^{-1}\Phi_{+}$. Define $B=( \beta_{\Delta
}^{-1},-I_{p})
, $ where $\beta_{\Delta}$ is the $p\times p$ diagonal matrix such
that $%
\beta_{\Delta}1_{p}=\beta$ and $1_{p}$ is the $p\times1$ vector of $
1$'s.\vspace*{-2pt}

\begin{pf*}{Proof of Theorem \protect\ref{th1}}
We have $\hat{\theta}$
$=\arg
\min_{h\in\Theta}R(h) $, $\hat{\beta}=\hat{\beta}(\hat{
\theta})$, where
\[
R(h)=Q(\hat{\beta}(h),h),\qquad \hat{\beta
}(h)=M(h,h)^{-1}
\{ M(h,\theta)\beta+m(h)\}
\]
for $M(g,h)=\sum_{u\in\mathbb{N}}f(u;g)f(u;h)^{\prime}$, $m(h)=\sum
_{u\in
\mathbb{N}}f(u;h)x_{u}$. The subsequent proof implies that after
suitable\vadjust{\goodbreak}
norming $M(h,h)$ is well conditioned for relevant $h$ and large~$N$.  In
Theorem \ref{tha}, take $U(h)=\beta^{\prime}D\Psi(h)D\beta$,
$V(h)=V_{1}(h)-%
\{ V_{2}(h)-V_{2}(\theta)\} -\{ V_{3}(h)-V_{3}(\theta
)
\} , $ for $V_{1}(h)=\beta^{\prime}\{ P(h)-D\Psi
(h)D
\} \beta$, $V_{2}(h)=2m(h)^{\prime}M(h,h)^{-1}M(h,\theta
)\beta$, $V_{3}(h)=m(h)^{\prime}M(h,h)^{-1}m(h)$, with $\Psi(h)=\Phi
(\theta,\theta)-\Phi(\theta,h)\Phi(h,h)^{-1}\Phi(h,\theta)$, $%
P(h)=M(\theta,\theta)-M(\theta,h)M(h,h)^{-1}M(h,\theta)$. Define,
for $j=1,\ldots,p_{i}$, $i=1,\ldots,d$, and a finite $W$, positive scalars
$C_{ijw}$, $w=1,\ldots,W$, such that $C_{ijw}\leq C_{ij,w+1}$ for each such $w$.  Define
%
\begin{equation}\label{eqB.1}
C_{w}=( C_{11w},\ldots,C_{1p_{1}w},\ldots,C_{d1w},\ldots,C_{dp_{d}w})
,\qquad%
 w=1,\ldots,W+1.\vspace*{-2pt}
\end{equation}
Define neighbourhoods $\mathcal{N}_{ij}(C_{ijw})=\{ h_{ij}\dvt
\vert
h_{ij}-\theta_{ij}\vert<C_{ijw}\} $, $j=1,\ldots,p_{i}$, $%
i=1,\ldots,d$, $w=1,\ldots,W+1$. Finally, define for $w=1,\ldots,W+1$,%
%
\begin{equation} \label{eqB.2}
\mathcal{N}(C_{w})=\prod_{i=1}^{d}\prod_{j=1}^{p_{i}}\mathcal{N}_{ij}(C_{ijw}),\vspace*{-2pt}
\end{equation}
and then $\bar{\mathcal{N}}(C_{w})$, $\mathcal{S}_{w}$ as in (\ref{eqA.1}).  Take $%
C_{ij1}=N^{\chi-\zeta_{ij}-{1/2}}\sim B_{i}^{\theta_{ij}}N^{\chi-{1/2}
}n_{i}^{-\theta_{ij}}$, $j=1,\ldots,\allowbreak p_{i}$, $i=1,\ldots,d$, so we need to show
that $P( \hat{\theta}\in\bar{\mathcal{N}}(C_{1}))
\rightarrow0$
as $N\rightarrow\infty$. We check~(i) and~(ii) of Theorem \ref{tha}, where
(\ref{eqA.2}) reduces to the requirement $\inf_{h\in\mathcal
{S}_{w}}U(h)/s_{w}>\eta
$ for large enough $N$  and $\eta$ as in (\ref{eqA.2}).  From (\ref{eqB.1}) and (\ref{eqB.2}),
\[
\mathcal{S}_{w}\subset\Theta\cap\mathcal{T}_{w},\vspace*{-2pt}
\]
where
\[
\mathcal{T}_{w}=\bigcup^{d}_{i=1}
\bigcup^{p_{i}}_{{j=1}}\{ h_{ij}\dvt \vert h_{ij}-\theta_{ij}\vert
\geq
C_{ijw};h_{kl}\dvt h_{kl}\in(-1/2,\infty),\mbox{all }(
k,l)
\neq( i,j) \} .\vspace*{-2pt}
\]
It follows from Proposition \ref{pr1} that
\[
\inf_{h\in\mathcal{S}_{w}}U(h)\geq\eta^{\ast}N\min_{i,j}\beta_{ij}^{2}\sum_{i=1}^{d}
\sum_{j=1}^{p_{i}}n_{i}^{2\theta_{ij}}C_{ijw}^{2}\geq\frac{\eta
}{p}%
\sum_{i=1}^{d}\sum_{j=1}^{p_{i}}
N^{1+2\zeta_{ij}}C_{ijw}^{2}.
\]
Thus (\ref{eqA.2}) is satisfied when
%
\begin{equation}\label{eqB.3}
\sum_{i=1}^{d}\sum_{j=1}^{p_{i}}%
N^{1+2\zeta_{ij}}C_{ijw}^{2}\geq ps_{w}.\vspace*{-2pt}
\end{equation}
Next, (\ref{eqA.3}) is implied if
%
\begin{eqnarray}\label{eqB.4}
\sup_{h\in S_{w}}\vert V_{1}(h)\vert&=&\mathrm{o}(s_{w}),
 \\[-2pt]
 \sup_{h\in S_{w}}\vert V_{2}(h)-V_{2}(\theta)
\vert
&=&\mathrm{o}_{p}(s_{w}), \label{eqB.5} \\[-2pt]
\sup_{h\in S_{w}}\vert V_{3}(h)\vert&=&\mathrm{o}_{p}(s_{w}),
\label{eqB.6}\vspace*{-2pt}
\end{eqnarray}
as $N\rightarrow\infty$. Note that in (\ref{eqB.5}) we are considering the
difference $V_{2}(h)-V_{2}(\theta)$ for $h$ suitably close to $\theta
$ and
this closeness\vadjust{\goodbreak} is important in obtaining the desired result, whereas in
the usual kind of consistency proof, for standard, non-mixed rate settings,
one more simply shows the convergence to zero in probability of a suitably
normalized $V_{2}(h)$, uniformly in $h\in$ $\Theta$. Now (\ref{eqB.6}) follows
from Proposition \ref{pr4}, while (\ref{eqB.4}) and (\ref{eqB.5}) follow from Propositions \ref{pr2}
and \ref{pr3},
respectively, if
%
\begin{equation}\label{eqB.7}
\sum_{i=1}^{d}\sum_{j=1}^{p_{i}}%
N^{1+2\zeta_{ij}-\delta^{\ast}}C_{ij,w+1}^{2}=\mathrm{o}(s_{w}),\vspace*{-2pt}
\end{equation}
where $\delta^{\ast}=\min[ \min_{1\leq i\leq d}\{
b_{i}/2+\min
( b_{i}\underaccent\bar{\Delta}_{i},0) \} ,2\chi] $
(implying $\delta^{\ast}>0)$ and
%
\begin{equation} \label{eqB.8}
\sum_{i=1}^{d}\sum_{j=1}^{p_{i}}N^{{1/2}
+\zeta_{ij}+\varepsilon}C_{ij,w+1}=\mathrm{o}(s_{w})\vspace*{-2pt}
\end{equation}
for some $\varepsilon>0$.

It remains to show\vspace*{1pt} that we can choose $W$ and the $s_{w}$, $C_{ijw}$, to
satisfy (i) of Theorem~\ref{tha} and (\ref{eqB.3}), (\ref{eqB.7}) and (\ref{eqB.8}). Now (\ref{eqB.3}) holds
for $w=1 $ if $s_{1}=N^{2\chi}$, and for $w>1$ if
\begin{eqnarray*}
s_{w} &=&s_{1}N^{(w-1)\delta^{\ast}/2}=N^{2\chi+(w-1)\delta^{\ast
}/2}, \\[-2pt]
C_{ijw} &=&C_{ij1}N^{(w-1)\delta^{\ast}/4}=N^{\chi-\zeta_{ij}-{1/2}
+(w-1)\delta^{\ast}/4},\qquad j=1,\ldots,p_{i},i=1,\ldots,d.\vspace*{-2pt}
\end{eqnarray*}
$ $Since
\[
N^{1+2\zeta_{ij}}C_{ij1}^{2}=s_{1},\qquad N^{1+2\zeta_{ij}-\delta^{\ast
}}C_{ij,w+1}^{2}=s_{1}N^{(w/2-1)\delta^{\ast}}=s_{w}N^{-\delta^{\ast}/2}\vspace*{-2pt}
\]
for all $i,j$, (\ref{eqB.7}) is satisfied. For all $i,j$,
\[
N^{{1/2}
+\zeta_{ij}+\varepsilon}C_{ij,w+1}=N^{\chi+\varepsilon+w\delta
^{\ast
}/4}=s_{w}N^{\varepsilon-\chi+\delta^{\ast}/4+(1-w)\delta^{\ast
}/4}=\mathrm{o}(s_{w}),\vspace*{-2pt}
\]
on taking $\varepsilon<\chi-\delta^{\ast}/4$, to satisfy (\ref{eqB.8}). Finally,
for all $i,j$, though $C_{ij1}\rightarrow0$ as $N\rightarrow\infty$ (no
matter how small $\delta^{\ast}$ or how large $\zeta_{ij})$, we have
$%
C_{ijw}\rightarrow\infty$ as $N\rightarrow\infty$ for large enough~$w$,
so there is a finite $W$ to satisfy (i) of Theorem \ref{tha}.\vspace*{-2pt}
\end{pf*}

\begin{pf*}{Proof of Theorem \protect\ref{th2}}  Omitted.\vspace*{-2pt}
\end{pf*}

\begin{pf*}{Proof of Theorem \protect\ref{th3}}
Put $a=(h^{\prime},b^{\prime})^{\prime}$, $Q(a)=Q(h,b)$ and define
$Q^{(1)}(a)=(\partial
/\partial a)Q(a)$, $Q^{(2)}(a)=(\partial/\partial a^{\prime})Q^{(1)}(a)$.
We have
\[
L_{+}Q^{(1)}(a)=-2\sum _{u\in\mathbb{N}}\{
y_{u}-b^{\prime
}f(u;h)\} H(u;h,b),\vspace*{-2pt}
\]
where $H(u;h,b)=[ ( L(u)f(u;h)\circ b) ^{\prime},(
Lf(u;h)) ^{\prime}] ^{\prime}$ with $L=$ $L(n)$ and $%
L_{+}Q^{(2)}(a)L_{+}=\sum _{i=1}^{3}Q_{i}^{(2)}(a)$, with%
\begin{eqnarray*}
Q_{1}^{(2)}(a) &=&2\sum _{u\in\mathbb
{N}}H(u;h,b)H(u;h,b)^{\prime},
\\[-2pt]
 Q_{2}^{(2)}(a) &=&2\sum _{u\in\mathbb{N}}\{ b^{\prime
}f(u;h)-\beta^{\prime}f(u;\theta)\} J(u;h,b),
\end{eqnarray*}
\begin{eqnarray*}
Q_{3}^{(2)}(a) &=&-2\sum _{u\in\mathbb{N}}x_{u}J(u;h,b),
\end{eqnarray*}
in which $J(u;h,b)$ is the 2$p\times2p$ $ $symmetric matrix with $(i,j)$th
$p\times p$ block $L(u)f_{\Delta}(u;\allowbreak h) L(u)b_{\Delta}$ for
$i=j=1$,  $%
L(u)f_{\Delta}(u;h)L$ for $i=1$, $j=2$ and $0$ for $i=j=2$, $b_{\Delta}$,
$f_{\Delta}(u,h)$ being the $p\times p$ diagonal matrices such that $%
b=b_{\Delta}1_{p}$, $f(u;h)=f_{\Delta}(u;h)1_{p}$.

By the mean value theorem
%
\begin{equation}\label{eqB.9}
D_{+}L_{+}^{-1}(\hat{\alpha}-\alpha)=\bigl( D_{+}^{-1}L_{+}\tilde{Q}%
^{(2)}L_{+}D_{+}^{-1}\bigr) ^{-1}D_{+}^{-1}L_{+}Q^{(1)}(\alpha),
\end{equation}
where $\tilde{Q}^{(2)}$ is formed from $Q^{(2)}(a)$ by evaluating its $i$th
row at $a=\bar{\alpha}_{(i)}$, where $\Vert\bar{\alpha
}_{(i)}-\alpha
\Vert\leq\Vert\hat{\alpha}-\alpha\Vert$, $%
i=1,\ldots,2p $.  By Proposition \ref{pr5} (\ref{eqB.9}) is
\[
\bigl\{ D_{+}^{-1}L_{+}Q^{(2)}(\alpha)L_{+}D_{+}^{-1}+\mathrm{O}_{p}(\log
N)^{-2}\bigr\} ^{-1}D_{+}^{-1}L_{+}Q^{(1)}(\alpha).
\]
Let $B_{\Delta}=\operatorname{diag}( \beta_{\Delta}^{-1},-I_{p}) $ and $%
\Gamma$ be the $2p\times2p$ matrix with $p\times p$ blocks $\Gamma
_{11}=0, $ $\Gamma_{21}=\Gamma_{12}^{\prime}=L^{-1}\Lambda$, $\Gamma
_{22}=-L^{-1}\Lambda-\Lambda^{\prime}L^{-1}$, with $\Lambda
=\Phi^{-1}\Phi_{+}\Upsilon^{-1}$. Noting Proposition \ref{pr6} and the
representations
\begin{eqnarray*}
BD_{+}^{-1}L_{+}Q^{(1)}(\alpha) &=&2N^{-
{1/2}}\sum _{u\in\mathbb{N}}\{ L(u)-L\}
D^{-1}f(u;\theta
)x_{u}, \\[-2pt]
B_{\Delta}\Gamma B_{\Delta}D_{+}^{-1}L_{+}Q^{(1)}(\alpha) &=&-2N^{-
{1/2}}\sum _{u\in\mathbb{N}}\bigl[ ( \beta_{\Delta
}^{-1}\Lambda
^{\prime}) ^{\prime},\bigl( L^{-1}\Lambda\{ L(u)-L
\}
-\Lambda^{\prime}\bigr) ^{\prime}\bigr] ^{\prime} \\[-2pt]
&&\hspace*{53pt}{}\times D^{-1}f(u;\theta)x_{u},
\end{eqnarray*}
we obtain from (\ref{eqB.9})
\begin{eqnarray*}
D_{+}L_{+}^{-1}(\hat{\alpha}-\alpha) &=&-N^{-
{1/2}
}B\sum _{u\in\mathbb{N}}[ \Upsilon^{-1}\{L(u)-L\}+\Lambda
^{\prime}] D^{-1}f(u;\theta)x_{u} \\[-2pt]
&&{}-N^{-
{1/2}
}\sum _{u\in\mathbb{N}}\bigl[ 0,\bigl( L^{-1}\Lambda
\{L(u)-L\}D^{-1}f(u;\theta)x_{u}\bigr) ^{\prime}\bigr] ^{\prime} \\[-2pt]
&&{}+\mathrm{O}_{p}( (\log N)^{-2}) N^{-{1/2}
}\sum _{u\in\mathbb{N}}( ( \beta_{\Delta}L(u))
^{\prime},L) ^{\prime}D^{-1}f(u;\theta)x_{u}.
\end{eqnarray*}
The last two terms are $\mathrm{O}_{p}((\log N)^{-1})$ by application of Lemmas \ref{le15}
and \ref{le10}, respectively. The proof is completed by applying
Proposition \ref{pr7} to the first term.\vspace*{-2pt}
\end{pf*}

\section{Propositions}\label{appc}\vspace*{-2pt}

\begin{proposition}\label{pr1}  For all $C_{w}$ given by \textup{(\ref{eqB.1})}
such that $C_{ijw}>0$, $j=1,\ldots,p_{i}$, $i=1,\ldots,d$, there exists
$\eta^{\ast}>0$ such that, for all $\theta\in\Theta$
\[
\inf_{h\in\bar{\mathcal{N}}(C_{w})}U(h)\geq\eta
^{\ast}N%
\sum_{i=1}^{d}\sum_{j=1}^{p_{i}}%
\beta_{ij}^{2}n_{i}^{2\theta_{ij}}C_{ijw}^{2}.\vadjust{\goodbreak}
\]
\end{proposition}

\begin{pf}
Non-singularity of $\Phi(h,h)$ for $h\in\Theta$, and
%
\begin{equation}\label{eqC.1}
\sup_{\Theta}\Vert\Phi(h,h)^{-1}\Vert\leq K,
\end{equation}
where $K$ throughout denotes a finite, positive generic constant, follow
from Lemmas \ref{le2} and \ref{le3}, numerators of elements of the inverse being
bounded and
denominators bounded away from zero. Now $\Psi(h)=[ (
I_{p},0) \Xi(h)^{-1}( I_{p},0) ^{\prime}] ^{-1}$,
where the $2p\times2p$ matrix $\Xi(h)$ has $(i,j)$th $p\times p$
submatrix $\Phi(\theta1(i=1)+h1(i=2),\theta1(i=1)+h1(i=2))$, $1(\cdot)$
denoting the indicator function and $\Xi(h)^{-1}$ existing on
$\bar{\mathcal{N}}(C_{w})$ as implied below. Introduce the $2p\times2p$
orthogonal permutation matrix $\Pi$ defined by $\Pi(1_{2}\otimes
a)=(
(1_{2}^{\prime}\otimes a_{1}^{\prime}),\ldots,(1_{2}^{\prime}\otimes
a_{d}^{\prime})) ^{\prime}$, for any $p\times1$ vector $a$
with $i$th $ p_{i}\times1$ subvector $a_{i}$. Then $\Pi\Xi(h)\Pi^{\prime}$
has the form of $T$ in Lemma \ref{le2} or \ref{le3}.

In the Lemma \ref{le2} situation, where no $\theta_{ij}$ is zero and no
$h_{ij}$ is
zero on $\bar{\mathcal{N}}(C_{w})$, we have $r_{i}=2p_{i}$,
$r=2p$, and
$v_{ik}=\theta_{ik}$, $k=1,\ldots,p_{i}$, $v_{ik}=\theta_{i,k-p_{i}}$, $%
k=p_{i}+1,\ldots,2p_{i}$. Denoting $E_{i}(h)=\operatorname{diag}\{ \theta
_{i1}-h_{i1},\ldots,\theta_{ip_{i}}-h_{ip_{i}}\} $ and $%
e_{i}(h)=\operatorname{diag}\{ E_{i}(h),-E_{i}(h)\} $, $e(h)=\operatorname{diag}\{
e_{1}(h),\ldots,e_{d}(h)\} $, inspection of the results of Lemma \ref{le2}
indicates that we may write $( \Pi\Xi(h)\Pi^{\prime})
^{-1}=e(h)^{-1}Ge(h)^{-1}$, where the $p\times p$ matrix $G$ is non-singular
and bounded on $\bar{\mathcal{N}}(C_{w})$. Then
\[
\Psi(h)=( I_{p},0) \Pi^{\prime}e(h)^{-1}Ge(h)^{-1}\Pi
(
I_{p},0) ^{\prime}=E(h)\tilde{G}^{-1}E(h),
\]
where $E(h)=\operatorname{diag}\{ E_{1}(h),\ldots,E_{d}(h)\} $, $\tilde{G}%
=( I_{p},0) \Pi^{\prime}G\Pi( I_{p},0)
^{\prime}$.
Thus $U(h)=\beta^{\prime}DE(h)\times \tilde{G}^{-1}E(h)D\beta
\geq
\beta^{\prime}D^{2}E(h)^{2}\beta/\operatorname{tr}(\tilde{G})$, whence the result
follows by boundedness of $\tilde{G}$ and $\inf_{h_{ij}\in\bar{%
\mathcal{N}}(C_{ijw})}( \theta_{ij}-h_{ij}) ^{2}=C_{ijw}^{2}$.

The details in the Lemma \ref{le3} setting, in which either $\theta_{ij}=0$
for one
$(i,j)$, or $h_{ij}$ can be zero on $\bar{\mathcal{N}}(C_{w})$ for
one $%
(i,j)$, are too similar to warrant inclusion.\vspace*{-2pt}
\end{pf}

\begin{proposition}\label{pr2}
%
\begin{equation} \label{eqC.2}
\sup_{h\in\mathcal{N}(C_{w})}\vert V_{1}(h)
\vert
\leq K\sum_{i=1}^{d}\sum_{j=1}^{p_{i}} \mathit{N}^{1+2\zeta_{ij}-\delta^{\ast
}}\mathit{C%
}_{ijw}^{2}.\vspace*{-2pt}
\end{equation}
\end{proposition}

\begin{pf}
Define $D(h)=N\operatorname{diag}\{
n_{1}^{h_{11}},\ldots,n_{1}^{h_{1p_{i}}},\ldots,n_{d}^{h_{d1}},\ldots,n_{d}^{h_{d}p_{d}}
\}
$, so $D=D(\theta)$, and $\tilde{M}(g, h)=D(g)^{-1}M(g,h)D(h)^{-1}$,
also $%
F_{1}(h)=\tilde{M}(\theta,\theta)-\tilde{M}(\theta,h)-\tilde
{M}(h,\theta
)+\tilde{M}(h,h)$, $F_{2}(h)=\{ \tilde{M}(\theta,h)-\tilde{M}%
(h,h)\} \tilde{M}(h,h)^{-1}\{ \tilde{M}(h,\theta)-\tilde
{M}%
(h,h)\} $, so we have the identity $%
D^{-1}P(h)D^{-1}=F_{1}(h)-F_{2}(h). $ Likewise, $\Psi(h)=\Psi
_{1}(h)-\Psi
_{2}(h)$, where
\begin{eqnarray*}
\Psi_{1}(h) &=&\Phi(\theta,\theta)-\Phi(\theta,h)-\Phi
(h,\theta)+\Phi(h,h), \\[-2pt]
\Psi_{2}(h) &=&\{ \Phi(\theta,h)-\Phi(h,h)\} \Phi
(h,h)^{-1}\{ \Phi(h,\theta)-\Phi(h,h)\} .
\end{eqnarray*}
Thus $V_{1}(h)=V_{11}(h)-V_{12}(h)$, where $V_{1i}(h)=\beta^{\prime
}D
\{ F_{i}(h)-\Psi_{i}(h)\} D\beta$,  $i=1,2$. Now~$V_{1i}(h)$ is bounded by
\begin{eqnarray*}
&&KN\sum_{i=1}^{d}\sum_{j=1}^{p_{i}}\sum_{\ell=1}^{p_{i}}n_{i}^{\theta_{ij}+\theta
_{i\ell}}\Biggl\vert\frac{1}{n_{i}}\sum^{n_{i}}_{u_{i}=1}v_{ij}(u_{i}/n_{i})v_{i\ell}(u_{i}/n_{i})-\int
 _{0}^{1}v_{ij}(x)v_{i\ell}(x)\,\mathrm{d}x\Biggr\vert\\
  &&\quad{}+KN\sum_{i=1}^{d}\sum_{j=1}^{p_{i}}
 \sum_{k=1}^{d}\sum_{\ell=1}^{p_{k}}n_{i}^{\theta_{ij}{}}n_{k}^{\theta_{k\ell}}\vert
\tilde{%
v}_{ij}\tilde{v}_{k\ell}-\bar{v}_{ij}\bar{v}_{k\ell
}
\vert,
\end{eqnarray*}
where $v_{ij}(x)=v( x;\theta_{ij},h_{ij}) $ with $v$ defined just before Lemma \ref{le6}, $\bar{v}%
_{ij}=\int _{0}^{1}v_{ij}(x)\,\mathrm{d}x$, $\tilde
{v}_{ij}=n_{i}^{-1}%
\sum _{u_{i}=1}^{n_{i}}v_{ij}(u_{i}/n_{i})$. $ $From Lemma \ref{le8} the
first modulus is bounded by
\[
K\vert h_{ij}-\theta_{ij}\vert\vert h_{i\ell
}-\theta
_{i\ell}\vert(\log n_{i})^{2}/n_{i}^{1+\min(2\underaccent\bar
{\Delta}%
_{i},0)}\leq KN^{\delta^{\ast}}
\]
because $\log n_{i}\leq\log N$, $n_{i}^{1+\min(2\underaccent\bar{\Delta}%
_{i},0)}=( B_{i}N^{b_{i}}) ^{1+\min(2\underaccent\bar{\Delta}%
_{i},0)}\geq N^{2\delta^{\ast}}/K$. The second modulus is bounded by
%
\begin{equation}\label{eqC.3}
\vert\tilde{v}_{ij}-\bar{v}_{ij}\vert\vert
\tilde{v}_{k\ell}\vert+\vert\tilde{v}_{k\ell}-%
\bar{v}_{k\ell}\vert\vert\bar{v}_{ij}
\vert.
\end{equation}
Now
\[
\vert\tilde{v}_{k\ell}\vert\leq\Biggl\{ \frac
{1}{n_{k}}%
\sum _{u_{k}=1}^{n_{k}}v_{k\ell}( u_{k}/n_{k})
^{2}\Biggr\} ^{{1/2}},\qquad \vert\bar{v}_{ij}\vert\leq\biggl\{ \int
 _{0}^{1}v_{ij}( x) ^{2}\,\mathrm{d}x\biggr\} ^{{1/2}
},
\]
so from Lemmas \ref{le6}, \ref{le7} and \ref{le8}, (\ref{eqC.3}) is bounded by
\[
K\biggl\{ \frac{(\log n_{i})^{2}}{n_{i}^{1+\min(\underaccent\bar{\Delta
}_{i},0)}}%
+\frac{(\log n_{k})^{2}}{n_{k}^{1+\min(\underaccent\bar{\Delta
}_{k},0)}}
\biggr\} \vert h_{ij}-\theta_{ij}\vert\vert h_{k\ell
}-\theta
_{k\ell}\vert,
\]
and the expression in braces is bounded by $N^{-\delta^{\ast}}$. Thus by
elementary inequalities, $\sup_{h\in\mathcal{N}(C_{w})} \vert
V_{1i}(h)\vert$ has the bound (\ref{eqC.2}). Next, $F_{2}(h)-\Psi_{2}(h)$
is
%
\begin{eqnarray} \label{eqC.4}
&&\{ \tilde{M}(\theta,h)-\tilde{M}(h,h)-\Phi(\theta,h)+\Phi
(h,h)\} \tilde{M}(h,h)^{-1}\{ \tilde{M}(h,\theta)-\tilde
{M}%
(h,h)\} \nonumber\\
&&\quad{}+\{ \Phi(\theta,h)-\Phi(h,h)\} \{ \tilde{M}%
(h,h)^{-1}-\Phi(h,h)^{-1}\} \{ \tilde{M}(h,\theta)-\tilde
{M}%
(h,h)\}\\
&&\quad{}+\{ \Phi(\theta,h)-\Phi(h,h)\} \Phi(h,h)^{-1}\{
\tilde{M}(h,\theta)-\tilde{M}(h,h)-\Phi(h,\theta)+\Phi(h,h)
\} .\nonumber
\end{eqnarray}
The final factor times $D\beta$ has a norm bounded by
%
\begin{eqnarray}\label{eqC.5}
&&KN^{{1/2}
}\sum_{i=1}^{d}\sum_{j=1}^{p_{i}}%
\Biggl\vert\sum_{\ell=1}^{p_{i}}n_{i}^{\theta
_{i\ell
}}\Biggl\{ \frac{1}{n_{i}}\sum _{u_{k}=1}^{n_{k}}v_{ij}(
u_{i}/n_{i}) ( u_{i}/n_{i}) ^{h_{i\ell}}-\int
 _{0}^{1}v_{ij}( x) x^{h_{i\ell}}\,\mathrm{d}x\Biggr\}
\Biggr\vert
\nonumber
\\[-8pt]
\\[-8pt]
\nonumber
&&\quad{} +KN^{
{1/2}
}\sum_{i=1}^{d}\sum_{j=1}^{p_{i}}%
\Biggl\vert\sum_{k=1}^{d}\sum_{\ell=1}^{p_{k}}
n_{k}^{\theta_{k\ell}}\Biggl\{ \tilde{v}_{ij}\frac
{1}{%
n_{k}}\sum _{u_{k}=1}^{n_{k}}( u_{k}/n_{k})
^{h_{k\ell}}-%
\bar{v}_{ij}\int _{0}^{1}x^{h_{k\ell}}\,\mathrm{d}x\Biggr\}
\Biggr\vert
.
\end{eqnarray}
The first term in braces is
\[
\frac{1}{n_{i}}\sum _{u_{i}=1}^{n_{i}}v\biggl( \frac
{u_{i}}{n_{i}}%
;\theta_{ij}+h_{i\ell},h_{ij}+h_{i\ell}\biggr) -\int
 _{0}^{1}v( x;\theta_{ij}+h_{i\ell},h_{ij}+h_{i\ell}
) \,\mathrm{d}x.
\]
By Lemma \ref{le8}, this is bounded by\vspace*{-2pt}
\[
K\vert\theta_{ij}+h_{i\ell}-h_{ij}-h_{i\ell}\vert
N^{-\delta^{\ast}}\leq K\vert\theta_{ij}-h_{ij}\vert
N^{-\delta^{\ast}}.\vspace*{-2pt}
\]
After rearrangement as before, and application also of Lemma \ref{le8}, the second
term in braces in (\ref{eqC.5}) has the same bound. Thus (\ref{eqC.5}) is bounded over
$h\in
\mathcal{N}(C_{w})$ by $K\times \sum _{i=1}^{d}\sum
 _{j=1}^{p_{i}}N^{
{1/2}+\zeta_{ij}-\delta^{\ast}}C_{ijw}$. On the other hand, using Lemma \ref{le6},\vspace*{-2pt}
%
\begin{equation}\label{eqC.6}
\Vert\beta^{\prime}D\{ \Phi(\theta,h)-\Phi(h,h)
\}
\Vert\leq K\sum _{i=1}^{d}\sum _{j=1}^{n_{i}}N^{
{1/2}
+\zeta_{ij}}C_{ijw}\vspace*{-3pt}
\end{equation}
uniformly in $h\in\mathcal{N}(C_{w})$. Using (\ref{eqC.1}), the contribution
to $%
V_{2i}(h)$ has the bound in (\ref{eqC.2}). To deal with the contributions from the
other two terms in (\ref{eqC.4}), standard manipulations indicate that it suffices
to show that\vspace*{-2pt}
%
\begin{eqnarray}
\sup_{h\in\mathcal{N}(C_{w})}\Vert\{ \tilde{M}%
(h,\theta)-\tilde{M}(h,h)\} D\beta\Vert&\leq&K\sum
 _{i=1}^{d}\sum _{j=1}^{n_{i}}N^{{1/2}
+\zeta_{ij}+\delta^{\ast}}C_{ijw}, \label{eqC.7} \\[-4pt]
\sup_{h\in\Theta}\Vert\tilde{M}(h,h)-\Phi
(h,h)
\Vert&\leq&KN^{-2\delta^{\ast}}. \label{eqC.8}\vspace*{-3pt}
\end{eqnarray}
Since the elements of $\tilde{M}(h,\theta)-\tilde{M}(h,h)$ are of form
$n_{i}^{-1}\sum _{u\in\mathbb{N}}( u_{i}/n_{i})
^{h_{ij}}v_{k\ell}( u_{k}/n_{k}) $, for $i=k$ or $i\neq k$,
(\ref{eqC.7}) follows much as before, using Lemmas \ref{le4} and \ref{le7}. Finally, (\ref{eqC.8}) is an
easy consequence of Lemma \ref{le5}.\vspace*{-3pt}
\end{pf}

\begin{proposition}\label{pr3}
For any $\varepsilon>0$\vspace*{-2pt}
%
\begin{equation}\label{eqC.9}
\sup_{h\in\mathcal{N}(C_{w})}\vert
V_{2}(h)-V_{2}(\theta
)\vert\leq K\sum _{i=1}^{d}\sum
_{j=1}^{n_{i}}N^{{1/2}
+\zeta_{ij}+\varepsilon}C_{ijw}.\vspace*{-3pt}
\end{equation}
\end{proposition}

\begin{pf}
We can write $V_{2}(h)-V_{2}(\theta)$ as\vspace*{-2pt}
%
\begin{eqnarray}\label{eqC.10}
&&2\{ m(h)-m(\theta)\} ^{\prime}\beta+2m(h)^{\prime
}M(h,h)^{-1}\{ M(h,\theta)-M(h,h)\} \beta \nonumber\\[-2pt]
&&\quad=2\{ \tilde{m}(h)-\tilde{m}(\theta)\} ^{\prime}D\beta\\[-2pt]
&&\qquad{}+2\tilde{m}(h)^{\prime}\tilde{M}(h,h)^{-1}\{ \tilde{M}(h,\theta)-%
\tilde{M}(h,h)\} D\beta, \nonumber\vspace*{-2pt}
\end{eqnarray}
where $\tilde{m}(h)=D(h)^{-1}m(h)$. Now\vspace*{-2pt}
%
\begin{equation} \label{eqC.11}
E\Bigl\{ \sup_{h\in\Theta}\Vert\tilde{m}(h)\Vert
\Bigr\}
\leq K\vspace*{-3pt}
\end{equation}
immediately from Lemma \ref{le10}. From the proof of Proposition \ref{pr2}, the last
term of
(\ref{eqC.10}) thus has the bound (\ref{eqC.9}). Next,\vspace*{-2pt}
\[
\vert\{ \tilde{m}(h)-\tilde{m}(\theta)\} ^{\prime
}D\beta\vert\leq K\sum _{i=1}^{d}\sum
 _{j=1}^{n_{i}}N^{\zeta_{ij}}\biggl\vert\sum _{u\in
\mathbb{N%
}}v_{ij}( u_{i}/n_{i}) x_{u}\biggr\vert\vspace*{-3pt}
\]
and by Lemma \ref{le11} its supremum over $\mathcal{N}(C_{w})$ has the bound
(\ref{eqC.9}).\vspace*{-2pt}\vadjust{\goodbreak}~%
\end{pf}

\begin{proposition}\label{pr4}\vspace*{-2pt}
\[
\sup_{h\in\Theta}\vert V_{3}(h)\vert\mathit
{\leq
K.}\vspace*{-2pt}
\]
\end{proposition}

\begin{pf}
Writing $V_{3}(h)=\tilde{m}(h)^{\prime}\tilde{M}(h,h)^{-1}\tilde{m}(h)$,
the result follows from (\ref{eqC.1}), (\ref{eqC.8}) and~(\ref{eqC.11}).\vspace*{-2pt}
\end{pf}

\begin{proposition}\label{pr5}As $N\rightarrow\infty$,\vspace*{-2pt}
\[
D_{+}^{-1}L\bigl\{ \tilde{Q}^{(2)}-Q^{(2)}(\alpha)\bigr\}
LD_{+}^{-1}=\mathrm{O}_{p}( (\log N)^{-2}) .\vspace*{-2pt}
\]
\end{proposition}

\begin{pf}
By elementary inequalities, the result follows if\vspace*{-2pt}%
\[
D_{+}^{-1}L\bigl\{ \tilde{Q}^{(2)}(\bar{\alpha})-Q^{(2)}(\alpha
)\bigr\}
LD_{+}^{-1}=\mathrm{O}_{p}( (\log N)^{-2})\vspace*{-2pt}
\]
for any $\bar{\alpha}$ such that $\Vert\bar{\alpha}-\alpha
\Vert\leq\Vert\hat{\alpha}-\alpha\Vert$. A typical
element of $Q_{1}^{(2)}(\bar{\alpha})-Q_{1}^{(2)}(\alpha)$ is\vspace*{-2pt}
%
\begin{eqnarray}\label{eqC.12}
&&2\sum_{u\in\mathbb{N}}(\log u_{i})^{\rho_{1}}(\log
n_{i})^{1-\rho_{1}}(\log u_{k})^{\rho_{2}}(\log n_{k})^{1-\rho_{2}}
\nonumber
\\[-10pt]
\\[-10pt]
\nonumber
&&\quad{}\times( \bar{\beta}_{ij}^{\rho_{1}}u_{i}^{\bar{\theta
}_{ij}}u_{k}^{%
\bar{\theta}_{k\ell}}\bar{\beta}_{k\ell}^{\rho_{2}}-\beta_{ij}^{\rho
_{1}}u_{i}^{\theta_{ij}}u_{k}^{\theta_{k\ell}}\beta_{k\ell}^{\rho
_{2}})\vspace*{-2pt}
\end{eqnarray}
for $i=k$ and $i\neq k$, and $\rho_{1},\rho_{2}=0,1$. We need to show
that (\ref{eqC.12}) $=\mathrm{O}_{p}( N^{1+\zeta_{ij}+\zeta_{k\ell}}/(\log
N)^{2}) $. With $\bar{\beta}_{ij},\bar{\beta}_{k\ell}$ replaced
by $%
\beta_{ij},\beta_{k\ell}$, it is bounded by\vspace*{-2pt}
%
\begin{equation}\label{eqC.13}
K(\log N)^{2}\frac{N}{n_{i}}\sum_{u_{i}=1}^{n_{i}}
\vert u_{i}^{\bar{\theta}_{ij}+\bar{\theta}_{i\ell}}-u_{i}^{\theta
_{ij}+\theta_{i\ell}}\vert,\qquad i=k,\vspace*{-2pt}
\end{equation}
or by\vspace*{-2pt}
%
\begin{equation}\label{eqC.14}
K(\log N)^{2}\frac{N}{n_{i}n_{k}}\sum_{u_{i}=1}^{n_{i}}%
\sum_{u_{i}=1}^{n_{i}}\vert u_{i}^{\bar
{\theta}%
_{ij}}u_{k}^{\bar{\theta}_{k\ell}}-u_{i}^{\theta_{ij}}u_{k}^{\theta
_{k\ell}}\vert,\qquad i\neq k.\vspace*{-2pt}
\end{equation}
Note that, for example,\vspace*{-2pt}
\[
\sum_{u_{i}=1}^{n_{i}}u_{i}^{\bar{\theta}%
_{ij}}=\mathrm{O}_{p}\Biggl( n_{i}^{\bar{\theta}_{ij}}\sup_{h_{ij}\in\Theta _{ij}}\Biggl\vert n_{i}^{-h_{ij}}
\sum_{u_{i}=1}^{n_{i}}%
u_{i}^{h_{ij}}\Biggr\vert\Biggr) =\mathrm{O}_{p}( n_{i}^{\theta
_{ij}}) ,\vspace*{-2pt}
\]
since $n_{i}^{\bar{\theta}_{ij}}=\mathrm{O}_{p}( n_{i}^{\theta_{ij}}\exp
(
N^{\chi-\zeta_{ij}-{1/2}}\log N) ) =\mathrm{O}_{p}( n_{i}^{\theta_{ij}}) $,
taking $\chi<\zeta_{11}+
{\frac12}$. Then from Theorem \ref{th2} and Lemma \ref{le12}, (\ref{eqC.13}) is\vspace*{-2pt}
\[
\mathrm{O}_{p}\bigl( (\log N)^{3}N^{1+\zeta_{ij}+\zeta_{i\ell}}( N^{\chi
-\zeta_{ij}-{1/2}
}+N^{\chi-\zeta_{i\ell}-
{1/2}}) \bigr) ,\vspace*{-2pt}
\]
which is $\mathrm{O}_{p}( N^{1+\zeta_{ij}+\zeta_{i\ell}}/(\log
N)^{2}) $
as desired, while using Lemma \ref{le4}, (\ref{eqC.14}) is\vspace*{-2pt}
\[
\mathrm{O}_{p}\bigl( (\log N)^{3}N^{1+\zeta_{ij}+\zeta_{k\ell}}( N^{\chi
-\zeta_{ij}-
{1/2}
}+N^{\chi-\zeta_{k\ell}-
{1/2}
}) \bigr) ,\vspace*{-2pt}
\]
which is $\mathrm{O}_{p}( N^{1+\zeta_{ij}+\zeta_{k\ell}}/(\log
N)^{2}) $
as desired.  Using Theorem \ref{th3}, it is readily seen that~(\ref{eqC.12})${}=
\mathrm{O}_{p}( N^{1+\zeta_{ij}+\zeta_{k\ell}}/(\log N)^{2})
$.\looseness=-1

The only elements of $Q_{2}^{(2)}(\bar{\alpha})-Q_{2}^{(2)}(\alpha)$ that
are not identically zero are the diagonal elements corresponding to the
three non-null submatrices in $J(u;h,b)$, and are of form
%
\begin{equation} \label{eqC.15}
2\sum_{u\in\mathbb{N}}\{ \bar{\beta}^{\prime
}f(u;\bar{%
\theta})-\beta^{\prime}f(u;\theta)\} u_{i}^{\bar{\theta}%
_{ij}}\{ (\log u_{i})^{2}\bar{\beta}_{ij}\} ^{\rho}(
\log
u_{i}\log n_{i}) ^{1-\rho}
\end{equation}
for $\rho=0,1$. We have to show this is $\mathrm{O}_{p}( N^{1+2\zeta
_{ij}}/(\log N)^{2}) $. After replacing $\bar{\beta}_{ij}$ by
$\beta
_{ij}$, it is bounded by
\[
K(\log N)^{2}\sum^{d}_{k=1}
\sum^{p_{k}}_{\ell=1}\biggl\vert\sum_{u\in\mathbb{N}}(u_{k}^{%
\bar{\theta}_{k\ell}}-u_{k}^{\theta_{k\ell}}) u_{i}^{\bar
{\theta}%
_{ij}}\biggr\vert.
\]
Proceeding much as before, this is
\begin{eqnarray*}
&&\mathrm{O}_{p}\Biggl( (\log N)^{2}\sum^{d}_{k=1}
\sum^{p_{k}}_{\ell=1} N^{1+\zeta_{k\ell}+\zeta_{ij}}N^{\chi
-\zeta
_{k\ell}-{1/2}
}\Biggr) \\[-2pt]
&&\quad=\mathrm{O}_{p}( (\log N)^{2}N^{
{1/2}
+\zeta_{ij}+\chi}) =\mathrm{O}_{p}\bigl( N^{1+2\zeta_{ij}}/(\log
N)^{2}\bigr) .
\end{eqnarray*}
Again, the same bound holds for (\ref{eqC.15}).

Finally, $Q_{3}^{(2)}(\bar{\alpha})-Q_{3}^{(2)}(\alpha)$ has non-zero
elements at the same locations, and they are of form
%
\begin{equation}\label{eqC.16}
-2(\log n_{i})^{1-\rho}\sum_{u\in\mathbb{N}}x_{u}\{
(%
\bar{\beta}_{ij}\log u_{i})^{\rho}u_{i}^{\bar{\theta}_{ij}}-(\beta
_{ij}\log u_{i})^{\rho}u_{i}^{\theta_{ij}}\}
\end{equation}
for $\rho=0,1$, which again will be shown to be $\mathrm{O}_{p}(
N^{1+2\zeta
_{ij}}/(\log N)^{2}) $. Replacing $\bar{\beta}_{ij}$ by $\beta_{ij}$
gives
\begin{eqnarray*}
&&-2\beta_{ij}(\log n_{i})^{1-\rho}\biggl\{ n_{i}^{\theta
_{ij}}\sum_{u\in\mathbb{N}}x_{u}(\log u_{i})^{\rho}v(
u_{i}/n_{i};\bar{%
\theta}_{ij},\theta_{ij}) \\[-2pt]
&&\hspace*{74pt}{} +( n^{\bar{\theta}_{ij}-\theta_{ij}}-1)
\sum_{u\in \mathbb{N}}x_{u}(\log u_{i})^{\rho}u_{i}^{\theta_{ij}}
\biggr\} \\[-2pt]
&&\quad =\mathrm{O}_{p}( (\log N)^{2}N^{{1/2}
+\zeta_{ij}+\chi}) =\mathrm{O}_{p}\bigl( N^{1+2\zeta_{ij}}/(\log
N)^{2}\bigr) ,
\end{eqnarray*}
applying Lemmas \ref{le10} and \ref{le11}, and $n^{\bar{\theta}_{ij}-\theta
_{ij}}-1=\mathrm{O}_{p}( (\log N)\vert\bar{\theta}_{ij}-\theta
_{ij}\vert) $.  We can show, as before, that (\ref{eqC.16}) has the
same bound.
\end{pf}

\begin{proposition}\label{pr6} As $N\rightarrow\infty$,
\[
D_{+}L^{-1}Q^{(2)}(\alpha)^{-1}L^{-1}D_{+}=
{\tfrac12}
B\Upsilon^{-1}B^{\prime}+
{\tfrac12}
B_{\Delta}\Gamma B_{\Delta}+\mathrm{O}_{p}( (\log N)^{-2}) .
\]
\end{proposition}

\begin{pf}
Clearly $Q_{2}^{(2)}(\alpha)\equiv0$. A typical non-zero element of $%
Q_{3}^{(2)}(\alpha)$ is
\[
-2\sum_{u\in\mathbb{N}}\{ (\log u_{i})^{2}\beta
_{ij}\} ^{\rho}( \log u_{i}\log n_{i}) ^{1-\rho
}u_{i}^{\theta_{ij}}x_{u}\vadjust{\goodbreak}
\]
for $\rho=0,1$, and from Lemma \ref{le10} this is
\[
\mathrm{O}_{p}( (\log N)^{2}N^{
{1/2}
+\zeta_{ij}}) =\mathrm{O}_{p}\bigl( N^{1+2\zeta_{ij}}/(\log N)^{2}\bigr)
\]
as desired. From Lemmas \ref{le13} and \ref{le14},%
\[
D_{+}^{-1}Q_{1}^{(2)}D_{+}^{-1}=2\operatorname{diag}\{ \beta_{\Delta
},I_{p}
\} \bigl( A+\mathrm{O}( (\log N)^{-2}) \bigr) \operatorname{diag}\{ \beta
_{\Delta
},I_{p}\} ,
\]
where $A$ has $p\times p$ submatrices $A_{ij}$ such that $A_{11}=L\Phi
L-L\Phi_{+}-\Phi_{+}^{\prime}L+\Phi_{++}$, $A_{12}=A_{21}^{\prime
}=L\Phi L-\Phi_{+}^{\prime}L$, $A_{22}= L\Phi L$. Thus $A^{-1}$
has $p\times p$ submatrices $A^{ij}$ such that $A^{11}=\Upsilon^{-1}$,
$%
A^{12}=A^{21^{\prime}}=\Lambda^{\prime}L^{-1}-\Upsilon^{-1}$, $%
A^{22}=L^{-1}\Phi^{-1}( L^{-1}+( \Phi L-\Phi_{+})
\Upsilon
^{-1}( L\Phi-\Phi_{+}) \Phi^{-1}L^{-1}) $.  It follows
that $A^{-1}=( I_{p},-I_{p}) \Upsilon^{-1}(
I_{p},-I_{p}) ^{\prime}+\Gamma$. The proof is straightforwardly
concluded.\looseness=-1
\end{pf}

\begin{proposition}\label{pr7} As $N\rightarrow\infty$,
\[
N^{-{1/2}
}\sum_{u\in\mathbb{N}}[ \Upsilon^{-1}\{
L(u)-L\} +\Lambda^{\prime}] D^{-1}f(u;\theta
)x_{u}\rightarrow_{d}\mathfrak{N}_{p}( 0,2\uppi F(0)\Upsilon
^{-1}) .
\]
\end{proposition}

\begin{pf}
Write $x_{u}=x_{u1}+x_{u2}$ for $x_{u1}=\sum _{v\in E_{M}}\xi
_{v}\varepsilon_{u-v}$, $x_{u2}=\sum _{v\in\bar{E}_{M}}\xi
_{v}\varepsilon_{u-v}$, with $E_{M}=\{ u\dvt \vert
u_{i}
\vert\leq M,i=1,\ldots,d\} $, $\bar{E}_{M}=\mathbb{Z}%
^{d}\setminus E_{M} $, for positive integer $M$. $ $For $\eta>0$,
choose $%
M $ such that $\sum_{v\in\bar{E}_{M}}\vert\xi_{v}
\vert
<\eta$. Writing
\[
g_{u}=[ \Upsilon^{-1}\{ L(u)-L\} +\Lambda^{\prime}%
] D^{-1}f(u;\theta),
\]
we have
%
\begin{eqnarray}\label{eqC.17}
E\biggl\Vert N^{-{1/2}
}\sum_{u\in\mathbb{N}}g_{u}x_{u2}\biggr\Vert^{2}
&=&N^{-1}%
\sum_{v\in\bar{E}_{M}}\sum_{w\in\bar{E}_{M}}\xi
_{v}\xi_{w}\sum_{u,u-v+w\in\mathbb{N}}g_{u}^{\prime}g_{u-v+w}
\nonumber
\\[-8pt]
\\[-8pt]
\nonumber
&\leq&\biggl( \sum_{v\in\bar{E}_{M}}\vert\xi
_{v}
\vert\biggr) ^{2}N^{-1}\sum_{u\in\mathbb{N}}\Vert
g_{u}\Vert^{2},
\end{eqnarray}
and
\[
\frac{1}{N}\sum_{u\in\mathbb{N}}\Vert g_{u}
\Vert
^{2}\leq\frac{K}{N}\sum_{u\in\mathbb{N}}
\sum^{d}_{i=1}
\sum^{p_{i}}_{j=1}[ 1+\{
\log
(u_{i}/n_{i})\} ^{2}] ( u_{i}/n_{i}) ^{2\theta
_{ij}}\leq K,
\]
by Lemmas \ref{le13} and \ref{le14}. Then (\ref{eqC.17}) $\leq K\eta^{2}$. Next write $N^{-{1/2}
}\sum _{u\in\mathbb{N}}g_{u}x_{u1}
=\break N^{-
{1/2}
}\sum _{w\in E^{\prime}}\varepsilon_{w}\times \sum _{u\in
E^{\prime\prime}}\xi_{u-w}g_{u}$, where
\begin{eqnarray*}
E^{\prime} &=&\{ w\dvt 1-M\leq w_{i}\leq n_{i}+M,
i=1,\ldots,d
\} , \\
E^{\prime\prime} &=&\{ u\dvt \max(1,w_{i}-M)\leq u_{i}\leq\min
(n_{i},w_{i}+M), i=1,\ldots,d\} .
\end{eqnarray*}
We may then apply a CLT, with $N$ and thus $N^{\ast}$ $=\prod
_{i=1}^{d}(n_{i}+2M)$ increasing, for independent random variables whose
squares are uniformly integrable. It remains to check two aspects. The
first is the Lindeberg condition,
\[
\frac{1}{N}\max_{w\in E^{\prime}}\biggl\Vert\sum_{u\in E^{\prime\prime}}\xi_{u-w}g_{u}\biggr\Vert^{2}\rightarrow0,
\qquad
\mbox{as }N\rightarrow\infty.
\]
The left side is bounded by
\[
\frac{K}{N}\max_{u\in\mathbb{N}}\Vert g_{u}
\Vert
^{2}\leq\frac{K}{N}\sum^{d}_{i=1}\sum^{p_{i}}_{j=1}\max_{u_{i}}[ \{ \log
(u_{i}/n_{i})\} ^{2}+1] ( u_{i}/n_{i}) ^{2\theta
_{ij}}\rightarrow0,
\]
since, for some $\eta>0$,
\[
( u_{i}/n_{i}) ^{2\theta_{ij}}\leq1(\theta_{ij}\geq
0)+n_{i}^{2\vert\theta_{ij}\vert}1( \theta
_{ij}<0) \leq N^{1-\eta},\qquad \vert\log
(u_{i}/n_{i})\vert\leq K\log N.
\]
The second aspect is the covariance structure:
%
\begin{equation}\label{eqC.18}
E\biggl\{ N^{-{1/2}}\sum_{u\in\mathbb{N}}g_{u}x_{1u}\biggr\} \biggl\{ N^{-{1/2}
}\sum_{u\in\mathbb{N}}g_{u}x_{1u}\biggr\} ^{\prime
}=N^{-1}%
\mathop{\sum\sum}_{v,w\in\bar{E}_{M}}\xi_{v}\xi_{w}\sum{}^{\prime}g_{u}g_{u+w-v},
\end{equation}
where the primed sum is over all $u$ such that $u$, $u+w-v\in\mathbb{N}$.
Since $M$ is fixed and $\Vert g_{u}\Vert\leq KN^{1-\eta}$,
for some $\eta>0$, (\ref{eqC.18}) differs by $\mathrm{o}(1)$, as $N\rightarrow\infty$,
from $N^{-1}\sum _{v,w\in E_{M}}\xi_{v}\xi_{w}\times \sum
 _{u\in\mathbb{N}}g_{u}g_{u+w-v}$. Using Lemma \ref{le16}, this
differs by $\mathrm{o}(1)$ from $N^{-1}( \sum _{v\in E_{M}}\xi
_{v}) ^{2}\sum _{u\in\mathbb{N}}g_{u}g_{u}^{\prime}$,
which, by Lemmas \ref{le13} and \ref{le14} and straightforward calculation and elimination,
equals
\begin{eqnarray*}
&&\biggl( \sum_{v\in E_{M}}\xi_{v}\biggr) ^{2}\{
\Upsilon
^{-1}\Phi_{++}\Upsilon^{-1}-\Upsilon^{-1}\Phi_{+}^{\prime}\Lambda
-\Lambda^{\prime}\Phi_{+}\Upsilon^{-1}+\Lambda^{\prime}\Phi
\Lambda
+\mathrm{O}(1/\log N)\} \\
&&\quad=\biggl( \sum_{v\in E_{M}}\xi_{v}\biggr) ^{2}\{
\Upsilon
^{-1}+\mathrm{O}(1/\log N)\} \rightarrow\biggl( \sum_{v\in E_{M}}%
\xi_{v}\biggr) ^{2}\Upsilon^{-1}
\end{eqnarray*}
as $N\rightarrow\infty$, and the last displayed expression differs by
$\mathrm{O}(\eta)$ from $( \sum _{\in\mathbb{Z}^{d}}\xi
_{v})
^{2}\Upsilon^{-1}=2\uppi F(0)\Upsilon^{-1}$.
\end{pf}

\section{Technical lemmas}\label{appd}
\begin{lemma}\label{le1}
Let $T$ be an $r\times r$ matrix, with $(i,j)$th $r_{i}\times r_{j}$ block
$T_{ij}$, $i,j=1,\ldots,d$, where $\sum
_{i=1}^{d}r_{i}=r$. Let $t_{i}$ be a column vector such that
$T_{ij}=t_{i}t_{j}^{\prime}$, $i\neq j$, and $T_{ii}-t_{i}t_{i}^{\prime}$
is positive definite,
$i,j=1,\ldots,d$. Then $T$ is non-singular, with $(i,j)$th
$r_{i}\times r_{j}$ submatrix
\begin{equation}\label{eqD.1}
T_{ii}^{-1}+\frac{T_{ii}^{-1}t_{i}t_{i}^{\prime}T_{ii}^{-1}}{1-\tau
_{i}}
\mathop{\sum_{s=1}^{d}}_{s\neq i}\frac{\tau_{s}}{
1-\tau_{s}}\Big/( 1+\sigma) ,\qquad i=j,
\end{equation}
and
%
\begin{equation} \label{eqD.2}
\frac{-T_{ii}^{-1}t_{i}t_{j}^{\prime}T_{jj}^{-1}}{(1-\tau_{i})(1-\tau
_{j})%
}\big/( 1+\sigma) ,\qquad i\neq j,
\end{equation}
where $\tau_{i}=t_{i}^{\prime}T_{ii}^{-1}t_{i}$,
$\sigma=\sum_{i=1}^{d}\tau_{i}/(1-\tau_{i})$.
\end{lemma}

\begin{pf}
Let $\tilde{T}$ be the $r\times r$ matrix with diagonal blocks $\tilde
{T}%
_{i}=T_{ii}-\tau_{i}\tau_{i}^{\prime}$, and zeros elsewhere, so
$T=\tilde{%
T}+tt^{\prime}$, where $t=(t_{1}^{\prime},\ldots,t_{d}^{\prime
})^{\prime}$.
Now because $\det\{T_{ii}-t_{i}t_{i}^{\prime}\}=\det\{T_{ii}\}(1-\tau
_{i})$, it follows that $\tau_{i}<1$, and
%
\begin{equation}\label{eqD.3}
\tilde{T}_{ii}^{-1}=T_{ii}^{-1}\{ I_{r_{i}}+(1-\tau
_{i})^{-1}t_{i}t_{i}^{\prime}T_{ii}^{-1}\} ,\qquad i=1,\ldots,d.
\end{equation}
Then $\tilde{T}^{-1}$ is the $r\times r$ matrix with diagonal blocks
$\tilde{%
T}_{ii}^{-1}$. Thus
%
\begin{equation}\label{eqD.4}
T^{-1}=\tilde{T}^{-1}\{ I_{r}-( 1+t^{\prime}\tilde{T}%
^{-1}t) ^{-1}tt^{\prime}\tilde{T}^{-1}\} .
\end{equation}
Now $t_{i}^{\prime}\tilde{T}_{ii}^{-1}=( 1+\tau_{i}( 1-\tau
_{i}) ^{-1}) t_{i}^{\prime}T_{ii}^{-1}=( 1-\tau
_{i}) ^{-1}t_{i}^{\prime}T_{ii}^{-1}$,  $i=1,\ldots,d$, and so $%
t^{\prime}\tilde{T}^{-1}=\{ ( 1-\tau_{1})
^{-1}t_{1}^{\prime}T_{11}^{-1},\ldots,( 1-\tau_{d})
^{-1}t_{d}^{\prime}T_{dd}^{-1}\} $, and thus $t^{\prime}\tilde
{T}%
^{-1}t=\sigma$. From (\ref{eqD.4}), the $(i,j)$th $r_{i}\times r_{j}$
submatrix of
$T^{-1}$, for $i\neq j$, is $-\tilde{T}_{ii}^{-1}t_{i}t_{j}^{\prime
}\tilde{T%
}_{jj}^{-1}/( 1+t^{\prime}\tilde{T}^{-1}t) $, which equals (\ref{eqD.2})
on substituting (\ref{eqD.3}), while for $i=j$ it is
\[
T_{ii}^{-1}+\frac{T_{ii}^{-1}t_{i}t_{i}^{\prime}T_{ii}^{-1}}{1-\tau
_{i}}-%
\frac{T_{ii}^{-1}t_{i}t_{i}^{\prime}T_{ii}^{-1}}{(1-\tau
_{i})^{2}}\big/
\{ 1+\sigma\} ,
\]
which equals (\ref{eqD.1}) after straightforward algebra.
\end{pf}

\begin{lemma}\label{le2}Let $T_{ii}$ be a
Cauchy matrix, having $(j,k)$th element
$(1+v_{ij}+v_{ik})^{-1}$, and let the $j$th element of
$t_{i}$ be $(1+v_{ij})^{-1}$, where $v_{ij}\in(-
{\frac12}
,\infty)\setminus\{0\}$, all $i,j$ and $v_{ij}\neq
v_{ik}$, for $j\neq k$. Then $T$ as defined in\vadjust{\goodbreak}
Lemma \ref{le1} is non-singular, and its inverse $T^{-1}$ has $(i,j)$th
$r_{i}\times r_{j}$ block\ with $(k,\ell)$th element
%
\begin{eqnarray}\label{eqD.5}
&&( 1+2v_{ik}) ( 1+2v_{i\ell}) \mathop{\prod^{r_{i}}_{m=1}}_{m\neq k}\frac{1+v_{ik}+v_{im}}{%
v_{ik}-v_{im}}\mathop{\prod^{r_{i}}_{m=1}}_{m\neq l}\frac{1+v_{i\ell}+v_{im}}{v_{i\ell}-v_{im}} \nonumber\\
 &&\quad{}\times\Biggl[ \frac{1}{1+v_{ik}+v_{i\ell}}-\biggl\{ \frac
{1}{v_{ik}(
1+v_{ik}) v_{i\ell}( 1+v_{i\ell}) }\prod^{r_{i}}_{m=1}\biggl( \frac{1+v_{im}}{v_{im}}\biggr)
^{2}\biggr\}
\\
\nonumber
&&\hspace*{8pt}\qquad\bigg/ \Biggl\{ \sum_{s=1}^{d}
\prod^{r_{i}}_{m=1}\biggl( \frac{1+v_{sm}}{v_{sm}}\biggr) ^{2}+1-d\Biggr\}
\Biggr] ,\qquad i=j,\\
\label{eqD.6}
&&\frac{( 1+2v_{ik}) ( 1+2v_{j\ell}) }{%
v_{ik}^{2}v_{j\ell}^{2}}\mathop{\prod^{r_{i}}_{m=1}}_{m\neq k}\frac{( 1+v_{ik}+v_{im}) ( 1+v_{im})
}{(
v_{im}-v_{ik}) v_{im}}
\nonumber
\\
&&\quad{}\times \mathop{\prod^{r_{i}}_{m=1}}_{m\neq l}\frac{( 1+v_{j\ell}+v_{jm}) (
1+v_{jm}) }{( v_{jm}-v_{j\ell}) v_{jm}}\\
&&\quad{}\bigg/
\Biggl\{ \sum_{s=1}^{d}
\prod^{r_{i}}_{m=1}\biggl( \frac{1+v_{sm}}{v_{sm}}\biggr) ^{2}+1-d\Biggr\} ,\qquad i\neq
j.\nonumber
\end{eqnarray}
\end{lemma}

\begin{pf}
From page 31 of Knuth \cite{Knu75}, $T_{ii}^{-1}$ has $(k,\ell)$th element
\begin{eqnarray*}
&&\prod^{r_{i}}_{m=1}( 1+v_{ik}+v_{im})
( 1+v_{i\ell}+v_{im}) \\
 &&\quad\bigg/\Biggl\{ ( 1+v_{ik}+v_{i\ell})
\mathop{\prod^{r_{i}}_{m=1}}_{m\neq l}( v_{ik}-v_{im})
\mathop{\prod^{r_{i}}_{m=1}}_{m\neq l}( v_{i\ell
}-v_{im}) \Biggr\} .
\end{eqnarray*}
For each $i$ define the $( r_{i}+1) \times(
r_{i}+1) $
non-singular Cauchy matrix $T_{ii}^{+}$ whose first $r_{i}$ rows are
$(T_{ii},t_{i})$ and whose last row is $( t_{i}^{\prime
},1)$. Thus, again from page~31 of Knuth \cite{Knu75}, the $(r_{i}+1,r_{i}+1)$th element
of its inverse is $(1-\tau_{i})^{-1}=\prod _{\ell
=1}^{r_{i}}( 1+v_{i\ell}^{-1}) ^{2}$. Thus
%
\begin{equation}\label{eqD.7}
1+\sigma=\prod^{r_{i}}_{\ell=1}( 1+v_{i\ell
}^{-1}) ^{2}+1-d.
\end{equation}
Also, the leading $r_{i}\times1$ subvector of the $(r_{i}+1)$th column
of $%
T_{ii}^{+-1}$ is $( 1-\tau_{i}) ^{-1}T_{ii}^{-1}t_{i}$, which
has $k$th element
\begin{eqnarray*}
&&( 1+v_{ik}) \prod^{r_{i}}_{m=1}(
1+v_{ik}+v_{im}) ( 1+v_{im})\bigg/\Biggl\{ (
1+v_{ik}) v_{ik}\mathop{\prod^{r_{i}}_{m=1}}_{m\neq l}( v_{ik}-v_{im}) \prod^{r_{i}}_{m=1}%
( -v_{im}) \Biggr\}
\\
&&\quad=\frac{1+2v_{ik}}{v_{ik}^{2}}\mathop{\prod^{r_{i}}_{m=1}}_{m\neq k}\frac{( 1+v_{ik}+v_{im}) ( 1+v_{im})
}{%
( v_{im}-v_{ik}) v_{im}}.
\end{eqnarray*}
The proof is completed by substitution and rearrangement.
\end{pf}

\begin{lemma}\label{le3}Let $T^{+}$ be the $(r+1)\times
(r+1)$
matrix whose first $r$ rows are $( T,t)
$ and whose last row is $( t^{\prime},1) $, with
$T$ and $t$ defined as in Lemmas \ref{le1} and \ref{le2}. Then
\[
T^{+-1}=\left[
\matrix{
T^{-1}\bigl( I_{r}+tt^{\prime}T^{-1}( 1-t^{\prime}T^{-1}t)
^{-1}\bigr) & -T^{-1}t( 1-t^{\prime}T^{-1}t) ^{-1} \vspace*{2pt}\cr
-( 1-t^{\prime}T^{-1}t) ^{-1}t^{\prime}T^{-1} & (
1-t^{\prime}T^{-1}t) ^{-1}%
}
\right] ,
\]
where $( 1-t^{\prime}T^{-1}t) ^{-1}=1+\sigma$.
\end{lemma}

\begin{pf}
From (\ref{eqD.1}) and (\ref{eqD.2})
\[
t^{\prime}T^{-1}t =\sum_{i=1}^{d}\biggl\{ \tau
_{i}+%
\frac{\tau_{i}^{2}}{(1+\sigma)(1-\tau_{i})}\biggl( \sigma-\frac{\tau
_{i}%
}{1-\tau_{i}}\biggr) \biggr\} -\frac{\sigma^{2}}{1+\sigma}+\frac{1}{(1+\sigma)}\sum
_{i=1}^{d}%
\frac{\tau_{i}^{2}}{(1-\tau_{i})^{2}}
=\frac{\sigma}{1+\sigma}
\]
after routine algebra. Thus $1-t^{\prime}T^{-1}t=1/(1+\sigma)$, and the
proof is readily completed.
\end{pf}

\begin{lemma}\label{le4}
For $\underaccent\bar{a}>-1$
\[
\sup_{a\geq\underaccent\bar{a}}\Biggl\vert\frac{1}{J}
\sum_{j=1}^{J}\biggl( \frac{j}{J}\biggr) ^{a}\Biggr\vert\leq K.
\]
\end{lemma}

\begin{pf}
The expression within the modulus is bounded by
\[
\int _{0}^{1}x^{a}\,\mathrm{d}x+1=\frac{1}{a+1}+1\leq\frac{1}{\underaccent\bar
{a}+1}%
+1\leq K.
\]
\upqed\end{pf}

\begin{lemma}\label{le5} For $\underaccent\bar{a}>-1$,
\[
\sup_{a\geq\underaccent\bar{a}}\Biggl\vert\frac{1}{J}
\sum_{j=1}^{J}\biggl( \frac{j}{J}\biggr) ^{a}-\frac{1}{1+a}
\Biggr\vert
\leq\frac{K}{J^{1+\min(\underaccent\bar{a},0)}}.
\]
\end{lemma}

\begin{pf}
The expression within the modulus is
%
\begin{equation}\label{eqD.8}
\frac{1}{J}\sum_{j=2}^{J-1}\int
 _{(j-1)/J}^{j/J}\biggl\{ \biggl( \frac{j}{J}\biggr)
^{a}-x^{a}
\biggr\} \,\mathrm{d}x+\frac{1}{J^{a+1}}-\int _{0}^{1/J}x^{a}\,\mathrm{d}x+\frac{1}{J}-\int
 _{1-1/J}^{1}x^{a}\,\mathrm{d}x.
\end{equation}
Using the mean value theorem, the first term in (\ref{eqD.8}) is bounded by
\[
\frac{2a}{J}\sum_{j=1}^{J}\biggl( \frac{j}{J}\biggr)
^{a-1}1( a>0) +\frac{\vert a\vert
}{J^{a+1}}\sum_{j=1}^{J}j^{a-1}1(a<0)\leq\frac{2}{J}+\frac{K}{J^{%
\underaccent\bar{a}+1}}.
\]
The last two integrals in (\ref{eqD.8}) are bounded by
\[
\frac{( 1/J) ^{a+1}}{a+1}+\frac{1}{a+1}\biggl\{ 1-\biggl(
1-\frac{1}{J}\biggr) ^{a}\biggr\} \leq\frac{K}{J^{\underaccent\bar{a}+1}}+\frac
{2}{J}.
\]
\upqed\end{pf}

Define, for $s\in\lbrack0,1]$, $v(s;a,b)=s^{a}-s^{b}$.

\begin{lemma}\label{le6}
For $\underaccent\bar{a}>-{\frac12}$
\[
\sup_{a,b\geq\underaccent\bar{a}}(a-b)^{-2}\int%
 _{0}^{1}v(x;a,b)^{2}\,\mathrm{d}x\leq K.
\]
\end{lemma}

\begin{pf}
The integral is
\[
\frac{1}{2a+1}-\frac{2}{a+b+1}+\frac{1}{2b+1}=\frac{2(a-b)^{2}}{%
(2a+1)(a+b+1)(2b+1)}\leq K(a-b)^{2}.
\]
\upqed\end{pf}

\begin{lemma}\label{le7}For $\underaccent\bar{a}>-
{\frac12}$,
%
\begin{equation} \label{eqD.9}
\sup_{a,b\in[\underaccent\bar{a},\bar{a}]}\Biggl\{
(a-b)^{-2}\sum_{j=1}^{J}v\biggl( \frac
{j}{J};a,b\biggr)
^{2}\Biggr\} \leq KJ(\log J)^{2}.
\end{equation}
\end{lemma}

\begin{pf}
By the mean value theorem,
%
\begin{equation}\label{eqD.10}
\vert v(s;a,b)\vert\leq s^{c}\vert\log s
\vert
\vert a-b\vert,\qquad s\in(0,1],
\end{equation}
where $\vert a-c\vert\leq\vert a-b\vert$. Also,
for such $c$,
%
\begin{equation}\label{eqD.11}
s^{c}\leq s^{\underaccent\bar{a}},\qquad s\in(0,1].
\end{equation}
Thus the quantity in braces in (\ref{eqD.9}) is bounded by
%
\begin{equation}\label{eqD.12}
K( \log J) ^{2}\sum_{j=1}^{J}\biggl(
\frac{j}{%
J}\biggr) ^{2\underaccent\bar{a}}\leq KJ(\log J)^{2},
\end{equation}
because $\underaccent\bar{a}>-{\frac12}$.
\end{pf}

\begin{lemma}\label{le8}For
$-1<\underaccent\bar{a}<\bar{a}<\infty$,
\[
\sup_{a,b\in[\underaccent\bar{a},\bar{a}]}\vert
a-b\vert^{-1}\Biggl\vert\frac{1}{J}\sum_{j=1}^{J}v\biggl( \frac{j}{J};a,b\biggr) -\int _{0}^{1}v(
x;a,b)
\,\mathrm{d}x\Biggr\vert\leq\frac{K(\log J)^{2}}{J^{1+\min(\underaccent\bar
{a},0)}}.
\]
\end{lemma}

\begin{pf}
The expression within the modulus is
%
\begin{equation}\label{eqD.13}
\sum_{j=2}^{J}\int _{(j-1)/J}^{j/J}\biggl\{
v\biggl( \frac{j}{J};a,b\biggr) -v( x;a,b) \biggr\} \,\mathrm{d}x+\frac
{1}{J}%
v\biggl( \frac{1}{J};a,b\biggr) -\int _{0}^{1/J}v(x;a,b)\,\mathrm{d}x.
\end{equation}
From (\ref{eqD.10}) and (\ref{eqD.11}), the last integral is bounded by
\[
K\int _{0}^{1/J}x^{\underaccent\bar{a}}\vert\log x\vert
\,\mathrm{d}x\vert a-b\vert\leq K(\log J)J^{-\underaccent\bar{a}-1}
\vert
a-b\vert,
\]
and the same bound results for the penultimate term of (\ref{eqD.13}). By the mean
value theorem $\vert v(s;a,b)-v(s-r;a,b)\vert$ is bounded by
%
\begin{equation}\label{eqD.14}
\vert s^{c}\log s-(s-r)^{c}\log(s-r)\vert\vert
a-b\vert,\qquad 0\leq r\leq1/J, s\geq2/J,
\end{equation}
where $\vert a-c\vert\leq\vert a-b\vert$, and
the first modulus is bounded by
%
\begin{eqnarray}\label{eqD.15}
&&\vert\{ s^{c}-(s-r)^{c}\} \log s\vert
+
\vert(s-r)^{c}\{ \log s-\log(s-r)\} \vert \nonumber
\\
&&\quad\leq s^{c}\vert\log s\vert\biggl\{ \biggl\vert1-
\biggl( 1-%
\frac{r}{s}\biggr) ^{c}\biggr\vert+\biggl( 1-\frac{r}{s}\biggr)
^{c}
\biggl\vert\log\biggl( 1-\frac{r}{s}\biggr) \biggr\vert\biggr\} \\
&&\quad\leq K\frac{s^{c-1}}{J}\vert\log s\vert.\nonumber
\end{eqnarray}
Thus the first term of (\ref{eqD.13}) is bounded by $\vert a-b
\vert$
times%
%
\begin{equation} \label{eqD.16}
\frac{K\log J}{J^{2}}\sum_{j=1}^{J}\biggl( \frac
{j}{J}%
\biggr) ^{\underaccent\bar{a}-1}
=\mathrm{O}\biggl( \frac{\log J}{J^{\underaccent\bar
{a}+1}}1(\underaccent\bar{a}\leq0)+\frac{(\log J)^{2}}{J}1(\underaccent\bar{a}=0)+\frac
{\log J%
}{J}1(\underaccent\bar{a}>0)\biggr) .\qquad
\end{equation}
\upqed\end{pf}

\begin{lemma}\label{le9}For $\underaccent\bar{a}>-{\frac12}$,
\begin{eqnarray*}
&&\mathop{\mathop{\sup_{t{a_{j},b_{i}\in[\underaccent\bar{a},\bar{a}]}}}_{i=1,2}}\Biggl\{ \prod^{2}_{i=1}\vert
a_{i}-b_{i}\vert\Biggr\} ^{-1}\Biggl\vert\frac{1}{J}\sum_{j=1}^{J}
\prod^{2}_{i=1}v\biggl( \frac
{j}{J}%
;a_{i},b_{i}\biggr)  -\int _{0}^{1}\prod^{2}_{i=1}
v(x;a_{i},b_{i})\,\mathrm{d}x\Biggr\vert\\
&&\quad\leq\frac{K(\log J)^{3}}{J^{1+\min(2%
\underaccent\bar{a},0)}}.
\end{eqnarray*}
\end{lemma}

\begin{pf}
The expression within the second modulus is
%
\begin{eqnarray}\label{eqD.17}
&&\sum_{j=2}^{J} \int
_{(j-1)/J}^{j/J}\Biggl\{
\prod^{2}_{i=1}v\biggl( \frac
{j}{J};a_{i},b_{i}\biggr) -%
\prod^{2}_{i=1}v(x;a_{i},b_{i})\Biggr\} \,\mathrm{d}x
\nonumber
\\[-8pt]
\\[-8pt]
\nonumber
&&\quad{}+J^{-1}\prod^{2}_{i=1}v\biggl( \frac{1}{J}%
;a_{i},b_{i}\biggr) -\int _{0}^{1/J}\prod^{2}_{i=1}
v(x;a_{i},b_{i})\,\mathrm{d}x.
\end{eqnarray}
Similarly to the proof of Lemma \ref{le7}, the last term is bounded by%
\[
K\int _{0}^{1/J}x^{2\underaccent\bar{a}}(\log x)^{2}\,\mathrm{d}x
\prod^{2}_{i=1}\vert a_{i}-b_{i}\vert\leq\frac{K(\log
J)^{2}}{J^{2\underaccent\bar{a}+1}}\prod^{2}_{i=1}\vert
a_{i}-b_{i}\vert.
\]
The expression in braces in (\ref{eqD.17}) can be written
\begin{eqnarray*}
&&\biggl\{ v\biggl( \frac{j}{J};a_{1},b_{1}\biggr) -v(
x;a_{1},b_{1}) \biggr\} v\biggl( \frac{j}{J};a_{2},b_{2}\biggr) \\
&&\quad{}+v( x;a_{1},b_{1}) \biggl\{ v\biggl( \frac{j}{J}%
;a_{2},b_{2}\biggr) -v( x;a_{2},b_{2}) \biggr\} .
\end{eqnarray*}
Both terms are treated similarly; we consider only the first. From the
bounds (\ref{eqD.14}), (\ref{eqD.15}) its first factor is bounded by
$(K/J)(j/J)^{\underaccent\bar{%
a}-1}(\log J)\vert a_{1}-b_{1}\vert$, and its second one
by $%
K(j/J)^{\underaccent\bar{a}}(\log J)\vert a_{2}-b_{2}\vert$. Thus
its contribution is
\[
\mathrm{O}\Biggl((\log J)^{2}J^{1+2\underaccent\bar{a}}\sum
_{j=1}^{J}j^{2\underaccent\bar{a}%
-1}\Biggr),
\]
whence the result follows by an analogous calculation to (\ref{eqD.16}).
\end{pf}

\begin{lemma}\label{le10}For $i=1,\ldots,d$ and $-
{\frac12}<\underaccent\bar{a}<\bar{a}<\infty$\textit{, and all} $q\geq0$
%
\begin{equation}\label{eqD.18}
E\biggl\{ \sup_{a\in[\underaccent\bar{a},\bar{a}]}
\biggl\vert N^{-
{1/2}
}\sum_{u\in\mathbb{N}}\biggl( \frac{u_{i}}{n_{i}}\biggr)
^{a}(\log u_{i})^{q}x_{u}\biggr\vert\biggr\} \leq K(\log N)^{q}.
\end{equation}
\end{lemma}

\begin{pf}
By summation by parts
\begin{eqnarray*}
&&\sum_{u_{i}=1}^{n_{i}}\biggl( \frac{u_{i}}{n_{i}}%
\biggr) ^{a}(\log u_{i})^{q}x_{u} \\
&&\quad=\sum_{u_{i}=1}^{n_{i}-1}\biggl\{ \biggl( \frac
{u_{i}}{%
n_{i}}\biggr) ^{a}-\biggl( \frac{u_{i}+1}{n_{i}}\biggr) ^{a}\biggr\}
\sum_{\ell=1}^{u_{i}}(\log\ell
)^{q}x_{u_{1},\ldots,\ell
,\ldots,u_{d}}+\sum_{u_{i}=1}^{n_{i}}(\log u_{i})^{q}x_{u},
\end{eqnarray*}
where $x_{u_{1},\ldots,\ell,\ldots,u_{d}}$ is $x_{u}$
with $%
u_{i}$ replaced by $\ell$. Thus the expression in the modulus in~(\ref{eqD.18}) is
%
\begin{equation}\label{eqD.19}
N^{-
{1/2}
}\sum_{u_{i}=1}^{n_{i}-1}\biggl( \frac{u_{i}}{n_{i}}
\biggr) ^{a}\{ 1-( 1+u_{i}^{-1}) ^{a}\}
H_{i}(u_{i})+n^{-
{1/2}}H_{i}(n_{i}),
\end{equation}
where
\[
H_{i}(s)=\mathop{\mathop{\sum^{n_{k}}_{u_{k}=1}}_{k=1,\ldots,d}}_{k\neq i}
\sum_{\ell=1}^{s}
x_{u_{1},\ldots,\ell,\ldots,u_{d}}(\log\ell)^{q}.
\]
The factor in braces in (\ref{eqD.19}) is bounded by $\vert a\vert
/u_{i}\leq K/u_{i}$, whereas $(u_{i}/n_{i})^{a}\leq(u_{i}/n_{i})^{%
\underaccent\bar{a}}$. Thus the left side of (\ref{eqD.18}) is bounded by
\begin{eqnarray*}
&&KN^{-{1/2}}\sum_{u_{i}=1}^{n_{i}-1}\biggl( \frac{u_{i}}{n_{i}}
\biggr) ^{\underaccent\bar{a}}\frac{1}{u_{i}}E\vert H_{i}(u_{i})
\vert
+n^{-
{1/2}
}E\vert H_{i}(n_{i})\vert\\
&&\quad \leq K(\log n_{i})^{q}n_{i}^{-
{1/2}
-\underaccent\bar{a}}\sum_{u_{i}=1}^{n_{i}}u_{i}^{\underaccent\bar{a%
}-
{1/2}
}+K(\log N)^{q}\leq K(\log N)^{q},
\end{eqnarray*}
since $\underaccent\bar{a}>-%
{\frac12}%
$ and\vspace*{-2pt}
\begin{eqnarray*}
EH_{i}(s)^{2} &=&\mathop{\mathop{\sum_{u_{k}=1}^{n_{k}}}_{k=1,\ldots,d}}_{k\neq i}
\sum_{\ell=1}^{s}
\mathop{\mathop{\sum_{v_{k}=1}^{n_{k}}}_{k=1,\ldots,d}}_{k\neq i}
\sum_{m=1}^{s}\gamma_{u_{1}-v_{1},\ldots,\ell
-m,\ldots,u_{d}-v_{d}}(\log\ell)^{q}(\log m)^{q} \\[-3pt]
&\leq&\frac{KNs}{n_{i}}(\log s)^{2q} \sum_{u\in\mathbb{Z}^{d}}\vert\gamma_{u}\vert\leq\frac{KNs(\log
s)^{2q}}{n_{i}%
}.\vspace*{-5pt}
\end{eqnarray*}
\upqed\end{pf}

\begin{lemma}\label{le11}
For $\underaccent\bar{a}>-{\frac12}$,\vspace*{-2pt}
%
\begin{equation} \label{eqD.20}
E\biggl\{ \sup_{a,b\in[\underaccent\bar{a},\bar{a}]}
\vert a-b\vert^{-1}\biggl\vert\sum_{u\in\mathbb{N}}
v( u_{i}/n_{i};a,b) x_{u}\biggr\vert\biggr\} \leq KN^{{1/2}}\log N.\vspace*{-2pt}
\end{equation}
\end{lemma}

\begin{pf}
By summation by parts,\vspace*{-2pt}
\[
\sum_{u_{i}=1}^{n_{i}}v(u_{i}/n_{i};a,b)x_{u}=
\sum_{u_{i}=1}^{n_{i}-1}\bigl\{ v(u_{i}/n_{i};a,b)-v\bigl((
u_{i}+1) /n_{i};a,b\bigr)\bigr\}
\sum _{\ell=1}^{u_{i}}x_{u_{1},\ldots,\ell,\ldots,u_{d}}.\vspace*{-2pt}
\]
From (\ref{eqD.14}) and (\ref{eqD.15}), the expression in braces is bounded by\vspace*{-2pt}
\[
K\biggl( \frac{\log n_{i}}{n_{i}}\biggr) \biggl( \frac
{u_{i}+1}{n_{i}}\biggr)
^{\underaccent\bar{a}-1}\leq\frac{K\log N}{u_{i}}\biggl( \frac{u_{i}}{n_{i}}%
\biggr) ^{\underaccent\bar{a}}.\vspace*{-2pt}
\]
Thus the left side of (\ref{eqD.20}) is bounded by
$K\log N\sum_{u_{i}=1}^{n_{i}-1}( u_{i}/n_{i}) ^{\underaccent\bar{a}%
}u_{i}^{-1}E\vert H_{i}(u_{i})\vert$, which, from the
proof of
Lemma \ref{le10} (with $q=0$), has the desired bound.\vspace*{-2pt}
\end{pf}

\begin{lemma}\label{le12}
Let $a>-{\frac12}$ be a scalar and $\tilde{a}=\tilde{a}_{J}$ be a sequence
such that $\tilde{a}-a=\mathrm{O}_{p}(J^{-\eta})$ as $J\rightarrow
\infty
$, for some $\eta>0$. Then for all $q\geq0$,\vspace*{-2pt}
\[
J^{-1-a}\sum_{j=1}^{J}(\log j)^{q}\vert
j^{\tilde{%
a}}-j^{a}\vert=\mathrm{O}_{p}( J^{-\eta}) ,\qquad \mbox{as }%
J\rightarrow\infty.\vspace*{-2pt}
\]
\end{lemma}

\begin{pf}
The left side is bounded by\vspace*{-2pt}
\begin{eqnarray*}
\sum_{j=1}^{J}(\log j)^{q}\biggl( \frac
{j}{J}\biggr)
^{a}\vert j^{\tilde{a}-a}-1\vert&\leq&\frac{1}{J}
\sum_{j=1}^{J}(\log j)^{q+1}\biggl( \frac{j}{J}\biggr) ^{a}
\vert
\tilde{a}-a\vert\\[-4pt]
 &\leq&KJ^{\eta/2}\mathrm{O}_{p}( J^{-\eta}) \frac{1}{J}
\sum_{j=1}^{J}\biggl( \frac{j}{J}\biggr) ^{a}=\mathrm{O}_{p}( J^{-\eta
/2}) .\vspace*{-6pt}
\end{eqnarray*}
\upqed\end{pf}

\begin{lemma}\label{le13}For $a>-{\frac12}$,
there is an $\eta>0$ such that for all sufficiently
large $J$,\vspace*{-2pt}
%
\begin{equation}\label{eqD.21}
\Biggl\vert\frac{1}{J}\sum_{j=1}^{J}( \log
j) \biggl( \frac{j}{J}\biggr) ^{a}-\frac{\log J}{a+1}+\frac
{1}{(a+1)^{2}%
}\Biggr\vert\leq KJ^{-\eta}.\vspace*{-2pt}\vadjust{\goodbreak}
\end{equation}
\end{lemma}

\begin{pf}
The left side is bounded by
%
\begin{equation}\label{eqD.22}
\frac{1}{J^{a+1}}\sum_{j=2}^{J}\int
 _{j-1}^{j}\vert(\log x)x^{a}-(\log j)j^{a}\vert
\,\mathrm{d}x+\biggl\vert\frac{1}{J^{a+1}}\int _{0}^{1}(\log
x)x^{a}\,\mathrm{d}x
\biggr\vert.
\end{equation}
The first modulus is bounded by
\begin{eqnarray*}
\vert\log x\vert\vert x^{a}-j^{a}\vert
+
\vert\log(x/j)\vert j^{a} &\leq&K(\log j)\{
(j-1)^{a-1}+j^{a-1}\} +j^{a-1} \\
 &\leq&K(\log j)j^{a-1}
\end{eqnarray*}
for $x\in\lbrack j-1,j]$, $j\geq2$. Thus the first term of (\ref{eqD.23}) is $
\mathrm{O}( (\log J)J^{-a-1}) $ for $a<0$, $\mathrm{O}( (\log
J)^{2}J^{-1}) $ for $a=0$, and $\mathrm{O}( (\log J)J^{-1}) $
for $%
a>0$. The last integral is $\mathrm{O}(J^{a-1})$. Since $a>-1$ there is an $\eta>0$
to satisfy (\ref{eqD.21}).
\end{pf}

\begin{lemma}\label{le14}
For any $a>-{\frac12}$, there is an
$\eta>0$ such that for all sufficiently
large $J$,
\[
\Biggl\vert\frac{1}{J}\sum_{j=1}^{J}(\log
j)^{2}\biggl(
\frac{j}{J}\biggr) ^{a}-\frac{(\log J)^{2}}{a+1}+\frac{2\log
J}{(a+1)^{2}}-%
\frac{2}{(a+1)^{3}}\Biggr\vert\leq J^{-\eta}.
\]
\end{lemma}

\begin{pf}
The left side is bounded by
\[
\frac{1}{J^{a+1}}\sum_{j=2}^{J}\int
 _{j-1}^{j}\vert(\log x)^{2}x^{a}-(\log j)^{2}j^{a}
\vert
\,\mathrm{d}x+\biggl\vert\frac{1}{J^{a-1}}\int _{0}^{1}(\log
x)^{2}x^{a}\,\mathrm{d}x\biggr\vert.
\]
The first integrand is bounded by
\[
(\log x)^{2}\vert x^{a}-j^{a}\vert+\vert\log
(x/j)\vert\vert\log( xj) \vert
j^{a}\leq
K(\log j)^{2}j^{a-1}
\]
as in the proof of Lemma \ref{le13}; the proof is completed in similar fashion.
\end{pf}

\begin{lemma}\label{le15}
For any $a>-{\frac12}$ and all sufficiently large $N$,
\[
E\biggl\{ N^{-{1/2}
}\sum_{u\in\mathbb{N}}\log( u_{i}/n_{i}) (
u_{i}/n_{i}) ^{a}x_{u}\biggr\} ^{2}\leq K.
\]
\end{lemma}

\begin{pf}
The left side is
\begin{eqnarray*}
&&N^{-1}%
\mathop{\sum\sum}_{u,v\in\mathbb{N}}
\biggl( \frac{u_{i}}{n_{i}}\biggr) ^{a}\biggl( \frac{v_{i}}{n_{i}}\biggr)
^{a}\log\biggl( \frac{u_{i}}{n_{i}}\biggr) \log\biggl( \frac
{v_{i}}{n_{i}}%
\biggr) \gamma_{u-v} \\
&&\quad\leq N^{-1}\sum_{u\in\mathbb{N}}\biggl( \frac
{u_{i}}{n_{i}}%
\biggr) ^{2a}\log^{2}\biggl( \frac{u_{i}}{n_{i}}\biggr) \sum_{v\in \mathbb{Z}^{d}}\vert\gamma_{u-v}\vert\leq K,
\end{eqnarray*}
by Assumption \ref{ass3} and straightforward application of Lemmas \ref{le13} and
\ref{le14}.
\end{pf}

\begin{lemma}\label{le16}For $a_{1},a_{2}>{\frac12}
$, $q_{1},q_{2}\geq0$, and any finite positive or negative
integer $M$, there is an $\eta>0$ such that for all sufficiently
large $J$,
%
\begin{eqnarray}\label{eqD.23}
&&\frac{1}{J}\sum_{j=1}^{J}\biggl\{ \log
^{q_{1}}\biggl(
\frac{j}{J}\biggr) \biggr\} \biggl( \frac{j}{J}\biggr) ^{a_{1}}\biggl\{
\log
^{q_{2}}\biggl( \frac{j+M}{J}\biggr) \biggl( \frac{j+M}{J}\biggr)
^{a_{2}}-\log^{q_{2}}\biggl( \frac{j}{J}\biggr) \biggl( \frac{j}{J}\biggr)
^{a_{2}}\biggr\}
\nonumber
\\[-8pt]
\\[-8pt]
\nonumber
&&\quad \leq\vert M\vert J^{-\eta}.
\end{eqnarray}
\end{lemma}

\begin{pf}
We have
\begin{eqnarray*}
\biggl\vert\biggl( \frac{j+M}{J}\biggr) ^{a_{2}}-\biggl( \frac
{j}{J}\biggr)
^{a_{2}}\biggr\vert&\leq&\frac{M}{j}\biggl( \frac{j}{J}\biggr)
^{a_{2}}, \\
 \biggl\vert\log^{q_{2}}\biggl( \frac{j+M}{J}\biggr) -\log^{q_{2}}\biggl(
\frac{j}{J}\biggr) \biggr\vert&\leq&\frac{M}{j}.
\end{eqnarray*}
By elementary inequalities the left side of (\ref{eqD.23}) is bounded by
\begin{eqnarray*}
&&\frac{KM(\log J)^{q_{1}+q_{2}}}{J^{a_{1}+a_{2}+1}}
\sum_{j=1}^{J}j^{a_{1}+a_{2}-1} \\
&&\quad\leq K\vert M\vert(\log J)^{q_{1}+q_{2}}\biggl\{ \frac{
1(a_{1}+a_{2}<0)}{J^{a_{1}+a_{2}+1}}+\frac{1(a_{1}+a_{2}=0)}{J}\log
J+\frac{
1(a_{1}+a_{2}>0)}{J}\biggr\} ,
\end{eqnarray*}
which is $\mathrm{O}(\vert M\vert J^{-\eta}).  $
\end{pf}
\end{appendix}

\section*{Acknowledgements}
I thank the Editor, Associate Editor and two referees for helpful comments.
I am also grateful to Francesca Rossi for carrying out the Monte Carlo
simulations. This research was supported by ESRC Grant RES-062-23-0036,
Spanish Plan Nacional de I+D+I Grant SEJ2007-62908/ECON, and a Catedra de
Excelencia at Universidad Carlos III, Madrid.


%

\printhistory


\begin{thebibliography}{25}

\bibitem{And71}
\begin{bbook}[mr]
\bauthor{\bsnm{Anderson},~\bfnm{T.~W.}\binits{T.W.}}
(\byear{1971}).
\btitle{The Statistical Analysis of Time Series}.
\baddress{New York}: \bpublisher{Wiley}.
\bid{mr={0283939}}
\end{bbook}
\endbibitem

\bibitem{Cre93}
\begin{bbook}[mr]
\bauthor{\bsnm{Cressie},~\bfnm{Noel A.~C.}\binits{N.A.C.}}
(\byear{1993}).
\btitle{Statistics for Spatial Data}.
\baddress{New York}: \bpublisher{Wiley}.
\bnote{Revised reprint of the 1991 edition}.
\bid{mr={1239641}}
\end{bbook}
\endbibitem

\bibitem{GaoLuTjs06}
\begin{barticle}[mr]
\bauthor{\bsnm{Gao},~\bfnm{Jiti}\binits{J.}},
  \bauthor{\bsnm{Lu},~\bfnm{Zudi}\binits{Z.}} \AND
  \bauthor{\bsnm{Tj{\o}stheim},~\bfnm{Dag}\binits{D.}}
(\byear{2006}).
\btitle{Estimation in semiparametric spatial regression}.
\bjournal{Ann. Statist.}
\bvolume{34}
\bpages{1395--1435}.
\bid{doi={10.1214/009053606000000317}, issn={0090-5364}, mr={2278362}}
\end{barticle}
\endbibitem

\bibitem{GirHidRob01}
\begin{barticle}[mr]
\bauthor{\bsnm{Giraitis},~\bfnm{L.}\binits{L.}},
  \bauthor{\bsnm{Hidalgo},~\bfnm{J.}\binits{J.}} \AND
  \bauthor{\bsnm{Robinson},~\bfnm{P.~M.}\binits{P.M.}}
(\byear{2001}).
\btitle{Gaussian estimation of parametric spectral density with unknown pole}.
\bjournal{Ann. Statist.}
\bvolume{29}
\bpages{987--1023}.
\bid{doi={10.1214/aos/1013699989}, issn={0090-5364}, mr={1869236}}
\end{barticle}
\endbibitem

\bibitem{GreRos84}
\begin{bbook}[mr]
\bauthor{\bsnm{Grenander},~\bfnm{Ulf}\binits{U.}} \AND
  \bauthor{\bsnm{Rosenblatt},~\bfnm{Murray}\binits{M.}}
(\byear{1984}).
\btitle{Statistical Analysis of Stationary Time Series}, \bedition{2nd} ed.
\baddress{New York}: \bpublisher{Chelsea}.
\bid{mr={0890514}}
\end{bbook}
\endbibitem

\bibitem{HalLuTra01}
\begin{barticle}[mr]
\bauthor{\bsnm{Hallin},~\bfnm{Marc}\binits{M.}},
  \bauthor{\bsnm{Lu},~\bfnm{Zudi}\binits{Z.}} \AND
  \bauthor{\bsnm{Tran},~\bfnm{Lanh~Tat}\binits{L.T.}}
(\byear{2001}).
\btitle{Density estimation for spatial linear processes}.
\bjournal{Bernoulli}
\bvolume{7}
\bpages{657--668}.
\bid{doi={10.2307/3318731}, issn={1350-7265}, mr={1849373}}
\end{barticle}
\endbibitem

\bibitem{HalLuYu09}
\begin{barticle}[mr]
\bauthor{\bsnm{Hallin},~\bfnm{Marc}\binits{M.}},
  \bauthor{\bsnm{Lu},~\bfnm{Zudi}\binits{Z.}} \AND
  \bauthor{\bsnm{Yu},~\bfnm{Keming}\binits{K.}}
(\byear{2009}).
\btitle{Local linear spatial quantile regression}.
\bjournal{Bernoulli}
\bvolume{15}
\bpages{659--686}.
\bid{doi={10.3150/08-BEJ168}, issn={1350-7265}, mr={2555194}}
\end{barticle}
\endbibitem

\bibitem{Jen69}
\begin{barticle}[mr]
\bauthor{\bsnm{Jennrich},~\bfnm{Robert~I.}\binits{R.I.}}
(\byear{1969}).
\btitle{Asymptotic properties of non-linear least squares estimators}.
\bjournal{Ann. Math. Statist.}
\bvolume{40}
\bpages{633--643}.
\bid{issn={0003-4851}, mr={0238419}}
\end{barticle}
\endbibitem

\bibitem{Jon43}
\begin{barticle}[auto:STB|2011/09/12|07:03:23]
\bauthor{\bsnm{Jones},~\bfnm{H.~T.}\binits{H.T.}}
(\byear{1943}).
\btitle{Fitting polynomial trends to seasonal data by the method of least
  squares}.
\bjournal{J. Amer. Statist. Assoc.}
\bvolume{38}
\bpages{453--465}.
\end{barticle}
\endbibitem

\bibitem{Knu75}
\begin{bbook}[vtex]
\bauthor{\bsnm{Knuth},~\bfnm{Donald~E.}\binits{D.E.}}
(\byear{1968}).
\btitle{The Art of Computer Programming. Vol. 1: Fundamental
Algorithms}.
\baddress{Reading, MA}: \bpublisher{Addison-Wesley}.
\end{bbook}
\endbibitem

\bibitem{Luetal07}
\begin{barticle}[mr]
\bauthor{\bsnm{Lu},~\bfnm{Zudi}\binits{Z.}},
  \bauthor{\bsnm{Lundervold},~\bfnm{Arvid}\binits{A.}},
  \bauthor{\bsnm{Tj{\o}stheim},~\bfnm{Dag}\binits{D.}} \AND
  \bauthor{\bsnm{Yao},~\bfnm{Qiwei}\binits{Q.}}
(\byear{2007}).
\btitle{Exploring spatial nonlinearity using additive approximation}.
\bjournal{Bernoulli}
\bvolume{13}
\bpages{447--472}.
\bid{doi={10.3150/07-BEJ5093}, issn={1350-7265}, mr={2331259}}
\end{barticle}
\endbibitem

\bibitem{Mal70}
\begin{barticle}[mr]
\bauthor{\bsnm{Malinvaud},~\bfnm{E.}\binits{E.}}
(\byear{1970}).
\btitle{The consistency of nonlinear regressions}.
\bjournal{Ann. Math. Statist.}
\bvolume{41}
\bpages{956--969}.
\bid{issn={0003-4851}, mr={0261754}}
\end{barticle}
\endbibitem

\bibitem{NagFul91}
\begin{barticle}[mr]
\bauthor{\bsnm{Nagaraj},~\bfnm{Neerchal~K.}\binits{N.K.}} \AND
  \bauthor{\bsnm{Fuller},~\bfnm{Wayne~A.}\binits{W.A.}}
(\byear{1991}).
\btitle{Estimation of the parameters of linear time series models subject to
  nonlinear restrictions}.
\bjournal{Ann. Statist.}
\bvolume{19}
\bpages{1143--1154}.
\bid{doi={10.1214/aos/1176348242}, issn={0090-5364}, mr={1126318}}
\end{barticle}
\endbibitem

\bibitem{Nie07}
\begin{barticle}[mr]
\bauthor{\bsnm{Nielsen},~\bfnm{Morten~{\O}rregaard}\binits{M.{\O}.}}
(\byear{2007}).
\btitle{Local Whittle analysis of stationary fractional cointegration and the
  implied-realized volatility relation}.
\bjournal{J. Bus. Econom. Statist.}
\bvolume{25}
\bpages{427--446}.
\bid{doi={10.1198/073500106000000314}, issn={0735-0015}, mr={2410027}}
\end{barticle}
\endbibitem

\bibitem{Rob08}
\begin{barticle}[mr]
\bauthor{\bsnm{Robinson},~\bfnm{P.~M.}\binits{P.M.}}
(\byear{2008}).
\btitle{Multiple local {W}hittle estimation in stationary systems}.
\bjournal{Ann. Statist.}
\bvolume{36}
\bpages{2508--2530}.
\bid{doi={10.1214/07-AOS545}, issn={0090-5364}, mr={2458196}}
\end{barticle}
\endbibitem

\bibitem{RobVid06}
\begin{barticle}[mr]
\bauthor{\bsnm{Robinson},~\bfnm{P.~M.}\binits{P.M.}} \AND
  \bauthor{\bsnm{Vidal~Sanz},~\bfnm{J.}\binits{J.}}
(\byear{2006}).
\btitle{Modified {W}hittle estimation of multilateral models on a~lattice}.
\bjournal{J. Multivariate Anal.}
\bvolume{97}
\bpages{1090--1120}.
\bid{doi={10.1016/j.jmva.2005.05.013}, issn={0047-259X}, mr={2276150}}
\end{barticle}
\endbibitem

\bibitem{SunPhi03}
\begin{barticle}[mr]
\bauthor{\bsnm{Sun},~\bfnm{Yixiao}\binits{Y.}} \AND
  \bauthor{\bsnm{Phillips},~\bfnm{Peter C.~B.}\binits{P.C.B.}}
(\byear{2003}).
\btitle{Nonlinear log-periodogram regression for perturbed fractional
  processes}.
\bjournal{J. Econometrics}
\bvolume{115}
\bpages{355--389}.
\bid{doi={10.1016/S0304-4076(03)00115-5}, issn={0304-4076}, mr={1984781}}
\end{barticle}
\endbibitem

\bibitem{Tjs78}
\begin{barticle}[mr]
\bauthor{\bsnm{Tj{\o}stheim},~\bfnm{Dag}\binits{D.}}
(\byear{1978}).
\btitle{Statistical spatial series modelling}.
\bjournal{Adv. in Appl. Probab.}
\bvolume{10}
\bpages{130--154}.
\bid{issn={0001-8678}, mr={0471224}}
\end{barticle}
\endbibitem

\bibitem{Tjs83}
\begin{barticle}[mr]
\bauthor{\bsnm{Tj{\o}stheim},~\bfnm{Dag}\binits{D.}}
(\byear{1983}).
\btitle{Statistical spatial series modelling. {II}. {S}ome further results on
  unilateral lattice processes}.
\bjournal{Adv. in Appl. Probab.}
\bvolume{15}
\bpages{562--584}.
\bid{doi={10.2307/1426619}, issn={0001-8678}, mr={0706617}}
\end{barticle}
\endbibitem

\bibitem{van00}
\begin{bbook}[mr]
\bauthor{\bparticle{van~de} \bsnm{Geer},~\bfnm{Sara~A.}\binits{S.A.}}
(\byear{2000}).
\btitle{Applications of Empirical Process Theory}.
\bseries{Cambridge Series in Statistical and Probabilistic Mathematics}
\bvolume{6}.
\baddress{Cambridge}: \bpublisher{Cambridge Univ. Press}.
\bid{mr={1739079}}
\end{bbook}
\endbibitem

\bibitem{Wu81}
\begin{barticle}[mr]
\bauthor{\bsnm{Wu},~\bfnm{Chien-Fu}\binits{C.F.}}
(\byear{1981}).
\btitle{Asymptotic theory of nonlinear least squares estimation}.
\bjournal{Ann. Statist.}
\bvolume{9}
\bpages{501--513}.
\bid{issn={0090-5364}, mr={0615427}}
\end{barticle}
\endbibitem

\bibitem{Yaj88}
\begin{barticle}[mr]
\bauthor{\bsnm{Yajima},~\bfnm{Yoshihiro}\binits{Y.}}
(\byear{1988}).
\btitle{On estimation of a regression model with long-memory stationary
  errors}.
\bjournal{Ann. Statist.}
\bvolume{16}
\bpages{791--807}.
\bid{doi={10.1214/aos/1176350837}, issn={0090-5364}, mr={0947579}}
\end{barticle}
\endbibitem

\bibitem{YajMat}
\begin{bmisc}[auto:STB|2011/09/12|07:03:23]
\bauthor{\bsnm{Yajima},~\bfnm{Y.}\binits{Y.}} \AND
  \bauthor{\bsnm{Matsuda},~\bfnm{Y.}\binits{Y.}}
  (\byear{2008}).
\bhowpublished{Asymptotic properties of the LSE of a spatial regression
  in both weakly and strongly dependent stationary random fields. Preprint
  CIRJE-F-587. Faculty of Economics, Univ. Tokyo}.
\end{bmisc}
\endbibitem

\bibitem{YaoBro06}
\begin{barticle}[auto:STB|2011/09/12|07:03:23]
\bauthor{\bsnm{Yao},~\bfnm{Q.}\binits{Q.}} \AND
  \bauthor{\bsnm{Brockwell},~\bfnm{P.~J.}\binits{P.J.}}
(\byear{2006}).
\btitle{Gaussian maximum likelihood estimation for ARMA models II: Spatial
  processes}.
\bjournal{Bernoulli}
\bvolume{12}
\bpages{403--429}.
\end{barticle}
\endbibitem

\end{thebibliography}
\end{document}